\documentclass[a4paper,12pt]{opelsarticle}
\usepackage{bm} %
\usepackage{lmodern}
\usepackage{xcolor}
\usepackage{lipsum}

\usepackage{amsmath,amsfonts,graphicx, pdfsync}
\usepackage[colorlinks=true, citecolor=red]{hyperref}
\usepackage{fullpage, graphicx}
\usepackage{tikzexternal}
\usepackage{filemod}
\DeclareGraphicsExtensions{.pdf, .jpg, .tif}

\definecolor{dbrown}{RGB}{0.41,0.05,0.05}
\definecolor{darkbrown}{RGB}{194,52,42}
\definecolor{darkred}{RGB}{250,30,30}

\def\n{\nabla}
\def\E{{\mathbb E}}
\def\p{\partial}
\def\ds{\displaystyle }

\def\RR{{\mathbb{R}}}

\def\N{{\mathbb{N}}}

\def\Q{{\mathcal Q}}
\def\R{{\mathbb R}}

\def\F{{\mathcal F}}

\def\P{{\mathcal P}}
\def\U{{\mathcal U}}
\def\V{{\mathcal V}}

\def\vd{{\mathbf d}}
\def\vr{{\mathbf r}}
\def\vu{{\mathbf u}}

\def\vz{{\mathbf z}}

\def\vF{{\mathbf F}}

\def\vX{{\mathbf X}}
\def\vY{{\mathbf Y}}
\def\vW{{\mathbf W}}
\def\vQ{{\mathbf Q}}

\def\va{{\mathbf a}}

\def\vp{{\mathbf p}}

\def\vD{{\mathbf D}}

\def\vd{{\mathbf d}}

\def\p{\partial}
\def\n{\nabla}

\def\d{\mathrm d}

\newtheorem{remark}{Remark}

\def\eq#1{\begin{eqnarray}#1\end{eqnarray}}
\def\eqs#1{\begin{eqnarray*}#1\end{eqnarray*}}

\begin{document}

\title{Performance  of a Markovian neural network versus dynamic programming on a fishing control problem}
\author{
Mathieu Lauri\`ere\footnote{\emph{lauriere@princeton.edu}, ORFE dept. Princeton University, USA },
Gilles Pag\`es\footnote{\emph{gilles.pages@sorbonne-universite.fr}, LPMA, Sorbonne Universit\'e, Paris, France },
Olivier Pironneau\footnote{\emph{olivier.pironneau@sorbonne-universite.fr }, LJLL, Sorbonne Universit\'e, Paris, France.}
}

\begin{frontmatter}
\begin{abstract}
Fishing quotas are unpleasant but efficient to control the productivity of a fishing site. A popular model has a stochastic differential equation for the biomass on which a stochastic dynamic programming or a Hamilton-Jacobi-Bellman algorithm  can be used to find the stochastic control -- the fishing quota. We compare the solutions obtained by dynamic programming against those obtained with a neural network which preserves the  Markov property of the solution.  The method is extended to a similar multi species model to check its robustness in high dimension.
 \end{abstract}
\end{frontmatter}

\section{Introduction}
Too much fishing can deplete the biomass to a disastrous level and leave fishermen out of work, or even, in some part of the world, starving. There are several ways to control fishing. One is to forbid fishing in regions and to alternate fishing and non-fishing zones \cite{moussaoui}. Another is by imposing quotas. In \cite{augerOP} statistical learning was shown to be very efficient to calibrate the parameters of the fishing model of \cite{MP18}. In \cite{pagespironneau} a stochastic control problem was derived from the model used in \cite{MP18} and to speed up the computation of optimal quotas, a solution by statistical learning was proposed and compared to standard  stochastic control solutions like the Hamilton-Jacobi-Bellman equations (HJB).  However the solution provided by the neural network was not Markovian as it used the future states, given by the stochastic model, to optimise the present.  In this article we propose to study the performance of a Markovian neural network, in the spirit of~\cite{GobetMunos05,HanE16,CarmonaLauriere19b}. 

The unpopularity of severe quotas is modeled by a penalty in the cost of an optimisation problem to compute the best fishing strategy which preserves the fish biomass, \textit{i.e.}, keeps it close to an ideal state $X_d$.
Quotas on fishing should also be as stable in time as possible because fishermen need to know that the quota will not move too much from one day to the next.  This is modeled by another penalty on the time variations of the quota.

Let $X_t$ be the fish biomass at time $t$, $E_t$ the fishing effort, interpreted as the number of boats at sea.  In \cite{MP18},  $q X_t$, with the {\it catchability} constant $q$, is the maximum weight of fish that a boat can mechanically catch,  meaning that the more fish there is the more fishermen will catch them, but proportionally to the capacity of his equipment.

Ideally a fisherman may want to catch a quantity $Q_t$ on day $t$. Imposing a quota $Q_M$ means $Q_t<Q_M$.  

\section{The fishing site model}

Consider a fishing site with $d$ types of fish in interaction, in the sense that some depend on others for food. Let $\vX(t)\in\R^d$ be the fish biomasses at time $t$, $E(t)$ the fishing effort -- interpreted as the number of boats at sea -- and the capacity to catch type $i$ is $q X_i(t)$, with the catchability $q$ constant over types for simplicity.

An optimal strategy with quotas  for each species $i$ is given to each fisherman, to impose the maximum weight of fish, $Q_i(t)\in\R^d$, of that species caught on a day $t$. Hence the total amount of fish of type $i$ caught on a day  is $E(t)\min(q X_i(t),Q_i(t)),~i=1,\dots,d$.  

The logistic equation for $\vX(t)$ says that the biomasses is a consequence of the natural growth or decay rate $\vr$, the long time limit $\underline{\bm\kappa}^{-1}\vr$ of $X$, where $\underline{\bm\kappa}$ is the $d\times d$ species interaction matrix, and the depletion due to fishing:
\begin{eqnarray}\label{CE1v}\ds
   \frac{\d\vX}{\d t}(t)= \vX(t) \cdot\underline{\bm\Lambda}[\vr  - \underline{\bm\kappa}\vX(t) - \min(q \vX(t),\vQ(t)) E(t)].
\end{eqnarray}
The operator $\vr\mapsto\underline{\bm\Lambda}[\vr]$ transforms a vector $\vr\in\R^d$ into a diagonal $d\times d$ matrix. For example, with $d=2$ and $\kappa_{12}<0, \kappa_{21}>0$, then, in (\ref{CE1v}), 
\[
\vX(t) \cdot\underline{\bm\Lambda}[\vr  - \underline{\bm\kappa}\vX(t)] = \big[ X_1(r_1-\kappa_{11}X_1+|\kappa_{12}|X_2), X_2(r_2-|\kappa_{21}|X_1-\kappa_{22}X_2) \big]^T
\]
which means that the first species lives on its own but profit from the second species because it eats it as shown by the equation for the second species which has a death rate augmented by $|\kappa_{21}|X_1>0$.  

Let $\vF(t)=\min(q \vX(t),\vQ(t))$.  The fishing effort $E$ is driven by the profit $\vp\cdot\vF$  minus the operating cost of a boat $c$, where $p_i$ is the price of fish of species $i$:
\eq{\label{CE2}&\ds
   \frac 1 E\frac{\d E}{\d t}=\vp\cdot\vF    - c. 
}
 
The price is driven by the difference between the demand $\vD(\vp)$ and the resource $\vF E$:
\eq{\label{CE3}&\ds
    \Phi\frac{\d\vp}{\d t}= \vD(\vp) - \vF E &\hbox{ with the demand } \vD(\vp)_i = \frac{a'_i}{1+\gamma p_i}
}
where $\Phi$ is the inverse time scale at which the fish market price adjusts.  When $\Phi << 1$,
(\ref{CE3}) may be approximated by:
\[
 \vD(\vp) - \vF E =0 ~\Rightarrow~ \gamma p_i F_i E =  a'_i-F_i  E ~\Rightarrow~ \vp\cdot\vF E=\frac1\gamma 
\mathrm{tr}[\va'-E\vF],
\]
where the trace operator is defined by $\mathrm{tr}[\va]:=\sum_{i=1}^d a_i$.
Let us denote $a=\sum_1^d a'_i/\gamma^2$, $\tilde q_i =  q_i\gamma$, $\tilde Q_i = Q_i\gamma$, and $\tilde E = E/\gamma$. Then the whole system \eqref{CE1v}--\eqref{CE2} rewrites:
\eq{
\label{eq:CE3_tmp}
	\frac{\d\vX}{\d t}= \vX  \cdot\underline{\bm\Lambda}[\vr  - \underline{\bm\kappa} \vX  - \min(\tilde q \vX , \tilde \vQ ) \tilde E] ,~~
\frac{\d \tilde E}{\d t}=  a-(\mathrm{tr}[\min(\tilde q \vX , \tilde \vQ )] +c)\tilde E,
}
where $\vX\cdot\bm\Lambda[\cdot]$ is matrix-vector product $\bm\Lambda[\cdot]\vX$. Since we are not going to use the original variables in the sequel, we drop the tildas and write $q_i, Q_i,$ and $E$ instead of $\tilde q_i, \tilde Q_i,$ and $\tilde E$.

Finally we use another change of variables to get rid of the $q$ variable: we replace $t$ by $t/q$, $\vQ$ by $q\vQ$ and $(\vr,\underline{\bm\kappa},\va, c)$ by  $q(\vr,\underline{\bm\kappa},\va, c)$. Then  the above system~\eqref{eq:CE3_tmp} is identical but now with $q=1$.  In the end,  with $a=\rm{tr}[\va]:=\sum_1^d a_i$, the whole system for the evolution of the fish biomass $\vX$ and the fishing effort $E$ is:
\eq{
\frac{\d\vX}{\d t}= \vX  \cdot\underline{\bm\Lambda}[\vr  - \underline{\bm\kappa} \vX   - \min(\vX ,\vQ ) E] ,~~
\frac{\d E}{\d t}=  a-(\rm{tr}[min(\vX ,\vQ )] +c) E.
}

\subsection{The Stochastic vector control problem}
To prevent fish extinction, a constraint is set on the total catch  $\min(q X_i(t),Q_i(t))E(t)),$ $i=1,\dots,d,$ per species.  The value of $\vQ(t)$ is found by solving an optimisation problem described below. It is expected that $Q_i(t)<qX_i(t)$ otherwise the policy of quota is  irrelevant in the sense that the fisherman is given a maximum allowed catch  which is  greater than what he could mechanically catch.

Let us add a constraint $\vQ(t)\le \vQ_M$ so that $\min(q \vX(t),\vQ(t))=\vQ(t)$. Denote
 $u_i(t)=E(t)Q_i(t)/X_i(t)$ the optimal fishing policy for species $i$ per unit mass. We formulate the optimal control problem in terms of $u_i$ rather than $Q_i$. Then $E$ no longer appears and the problem becomes:
\eq{\label{opt2}
\ds \min_{\vu\in \U} \left\{\bar J:=\int_0^T|\vX(t)-\vX_d(t)|^2\d t~:~  \frac{\d\vX}{\d t}= \vX \cdot\underline{\bm\Lambda}\left[\vr-\vu - \underline{\bm\kappa} \vX \right],~~\vX(0)=\vX^0 \right\}
}
where the control is in the constraint space
$$
\U=\{u_m\le u_i(t) \le u_M, ~~i=1,\dots, d, t \in [0,T]\},
$$
which reflects a desire to ensure a minimal fishing $u_m$ and a maximum one $u_M$ so as to guarantee that $u_i<qE(t)$ at all time.

\subsection{More constraints on quotas}

To avoid small fishing quotas we use penalty and add to the criteria $\ds-\int_0^T\sum_1^d\alpha_i u_i(t)$.
Moreover,  to avoid too many daily changes we penalise the quadratic variation of $\vu$ over the time period $(0,T)$, \textit{i.e.}, to $\bar J$  we add $\beta\cdot\E[\vu]_t^{0,T}$, where the quadratic variation is defined as:
 \[
	 [\vu]^{0,T}_t = \lim_{\|P\|} \sum_1^K|\vu_{t_k}-\vu_{t_{k-1}}|^2
\]
where $P$ ranges over partitions of the interval $(0,T)=\cup_k(t_{k-1},t_k)$ and the limit is in probability when $\max_k|t_k-t_{k-1}|\to 0$.
Here It\^o calculus \cite{BIC} tells us that:
 \[
 \E[\vu_i]^{0,T}_t = \sigma^2\int_0^T \E[|X_t\cdot\n_X \vu_i|^2] \d t.
 \]
Hence, an optimal policy in the presence of noise is a solution of 
\eq{\label{reform}&
\min_{\vu\in \U}\Big\{& \bar J:=\E\left[\int_0^T\left[|\vX(t)-\vX_d(t)|^2\d t  - \bm\alpha \cdot \vu + \bm\beta\cdot[\vu]_t^{0,T}\right] \d t \right] ~: 
\cr&&
\label{eq:OC-dynX}
\d\vX_t=\vX \cdot\underline{\bm\Lambda}\left[(\vr-\vu  - \underline{\bm\kappa}\vX)\d t +\underline{\bm\sigma}\d\vW_t \right],~~\vX(0)=\vX^0 \Big\}.
}
For the sake of clarity we assume $\beta_i=\beta$ for all $i$. %
\begin{remark}\label{remm1b}
Let $\vY=\underline{\bm\kappa}\vX$. A multiplication of the SDE~\eqref{eq:OC-dynX} by $\underline{\bm\kappa}$ leads to
\[
\d\vY_t=\vY_t \cdot \underline{\bm\kappa}\,\underline{\bm\Lambda}\left[(\vr-\vu  - \vY_t)\d t +\underline{\bm\kappa}\otimes\underline{\bm\sigma}\d\vW_t \right]\underline{\bm\kappa}^{-T}.
\]
 Thus, if all the terms on the diagonal of $\underline{\bm\Lambda}\left[(\vr-\vu  - \vY_t)\d t +\underline{\bm\kappa}\otimes\underline{\bm\sigma}\d\vW_t \right]$ are equal and if the initials conditions $\underline{\bm\kappa}\vX(0)$ are equal,  then the components of $\vY$ are independent and equal because $\underline{\bm\kappa}$ and $\underline{\bm\Lambda}$ commute. It happens only if $\underline{\bm\kappa}\otimes\underline{\bm\sigma}\d\vW_t={\bf 1}\d W_t$.
\end{remark} 

\subsection{Existence of solution}

First, the solution of the SDE  exists and is positive.  For clarity the proof is done for $d=1$:

Note that $X\mapsto (r-\kappa X - u_t)X$ is locally Lipschitz, uniformly in $t$ since $u_t\in[u_m,u_M]$. Hence for every realization  $W_t(\omega)$  there is a unique strong solution until a blow-up time $\tau$ which is a stopping time for the filtration $\F_t^w=\underline{\sigma}(W_s, s\le t, {\cal N}_s)$.  
 
 On $[0,\tau[$ by It\^o calculus,
\eqs{&&
 X_t = X^0 \exp\left(({r-\frac{\sigma^2}2)t - \int_0^t(\kappa X_s + u_s) \d s + \sigma W_t}\right),
}
so:
$$
	X_t\le X^0 e^{S_t(u_m)},\forall t\in(0,\tau],~\hbox{ where } S_t(v):=(r-\frac{\sigma^2}2)t  - v t + \sigma W_t.
$$
Hence  $\ds \int_0^\tau X_s\d s =+\infty$ is impossible unless  $\tau=+\infty$, $\P$ a.s. Therefore 
\eqs{&\ds 
X^0 e^{ S_t(u_M)} \exp\left(-\kappa\int_0^t e^{S_s(u_M)}\d s\right)
& \le X_t\le
 X^0 e^{S_t(u_m)}.
}
 
Let $y=\log x$ and {$v=r-\frac{\sigma^2}2-\kappa e^y - u(e^y,t)$}.
From \cite{lebris}, if 
$v\in W^{1,1}_{loc}(\R), ~ v/(1+|y|)\in L^1(\R)\cap L^\infty(\R),~ 
\p_y v\in L^\infty(\R)$, then
 $Y_t$ has a PDF $\rho\in L^\infty\left(L^2(\R)\cap L^\infty(\R)\right)\cap L^2\left(H^1(\R)\right)$
  and the following \emph{equivalent} control problem has a solution:
\eqs{&& \label{control5}
\min_{v\in\V} J(v):=\int_{-R}^R\d y\int_0^T\left[(e^y-X_d)^2 -\alpha u(y,t) + \beta |\p_t v(y,t)|^2\right]\rho(y,t)\d t~:
\\ &&
\p_t \rho  +\p_y (v\rho) - \p_{yy}[\frac{\sigma^2}2\rho] = 0,~\rho(y,0)=\rho^0(y),~\forall y\in\R,~\forall t\in]0,T[.
}
where $
\V=\{v:~v_m\leq v\leq v_M,~ \|\p_t v\|_{L^2(\Q)} \leq K\}
$.

\begin{remark}
 This result also gives a computational method, by discretizing the above and solving it as an optimisation problem.  Of the four methods discussed in this article, this is the most expensive. It was tested in \cite{MLOP} and shown not to be superior in precision to other methods.
\end{remark}

\section{Time discretization}
We consider a uniform grid with $M+1$ points in time. Let ${h}= \frac T M$, $t^m=m{h}$ and, for any $f$ defined on $[0,T]$, $f^m$ an approximation of $f(t^m)$. Define the Euler scheme:
\eq{\label{euler1}
\vX^{m+1}=\vX^m\left(1+\underline{\bm\Lambda}\left[(\vr-\vu^m  - \underline{\bm\kappa}\vX^m)h +\underline{\bm\sigma}\delta\vW^m\right]\right),
}
where $\delta W_i^m = W_i^{(m+1){h}}-W_i^{m{h}}\sim\sqrt{h}\N_{0,1}$, $i=1,\dots, d$.
Note that positivity may not be preserved by this scheme, but in practice negative values do not seem to appear.

A Monte-Carlo method is used with $K$ sample solutions of \eqref{euler1} to compute the cost:
\eqs{%
	J_K(\vu):= \sum_{m=1}^M\frac hK\sum_{k=1}^K \left[|\vX^m_k-\vX_d|^2 -\bm\alpha\cdot\vu_k^m + \frac{\beta }{h}|\vu_k^{m+1}-\vu_k^{m}|^2\right],
}
with the convention that $\vu_k^{M+1}=\vu_k^{M}$.
\begin{remark}\label{remm1}
Let $\vX\mapsto \vu^{m+1}(\vX)$ and $\vX\mapsto \vu^m(\vX)$ be two given differentiable functions. Let $\sigma\in\R^+$. The following holds ({\rm reading hint: $u(X+a)$ is $u$ at $X+a$, not $u$ times $(X+a)$}):

\begin{align*}
\vu^{m+1}(\vX^{m+1}) 
&= \vu^{m+1}\left(\vX^m\left(1+{h}(r-\kappa \vX^m -\vu^m) + \sigma\delta W^m\right)\right)
\\
&\approx \vu^{m+1}\left(\vX^m\left(1+{h}(r-\kappa \vX^m -\vu^m) \right)\right) + \vu'^{m+1}\sigma \vX^{m+1}\delta W^m, 
\end{align*}
where the derivative $\vu'^{m+1}$ is evaluated at $\vX^m$ or  at $\vX^m\left(1+{h}(r-\kappa \vX^m -\vu^m) + \sigma\delta W^m\right)$.  
Hence:
\begin{align*}
\E[|\vu^{m+1}-\vu^{m}|^2|\vX^m] 
&\approx \left|\vu^{m+1}\left(\vX^m\left(1+{h}(r-\kappa \vX^m -\vu^m) \right)\right)-\vu^m\right|^2 
\\
&~~+ h{\sigma^2}|\vX^{m+1}\cdot\n\vu'^{m+1}|^2,
\end{align*}
so:
\begin{align}
&\frac\beta h\E[\sum_{m=1}^M h|\vu^{m+1}-\vu^{m}|^2|\vX^m] \notag
\\
&
\approx
\frac\beta h\sum_{m=1}^M h \left[\left|
\vu^{m+1}\left(\vX^m\left(1+{h}(r-\kappa \vX^m -\vu^m)\right)\right)-\vu^m\right|^2 + h{\sigma^2}|\vX^{m+1}\cdot\n\vu^{m+1}|^2\right].
\label{qv}
\end{align}
Even though the first term in (\ref{qv}) is dominated by the second term, we may keep it for numerical convenience.
\end{remark}

\section{Stochastic Dynamic Programming (SDP)}
Consider the {\it value function} 
\eq{\label{dmj}&
V(\vX,t) = &\min_{\vu\in \U}\Big\{\E\left[\int_t^T\left[|\vX(t)-\vX_d(t)|^2\d t  - \bm\alpha \cdot \vu + \beta[\vu]_t^{t,T}\right] \d t \right]  ~:~
\cr&&
\d\vX_t = \vX\cdot\underline{\bm\Lambda}\left[(\vr-\vu  - \underline{\bm\kappa}\vX)\d t +\underline{\bm\sigma}\d\vW_t \right],~~\vX(t)=\vX\Big\}.
}
Let $\bm\zeta_{{h}}(\vX,\vu^m,\vz)$ be the next iterate of a numerical scheme  for the SDE starting at $\vX$. For instance:
\[
\bm\zeta_h(\vX,\vu,\vz) =\vX(1+\underline{\bm\Lambda}\left[(\vr-\vu  - \underline{\bm\kappa}\vX)h + \underline{\bm\sigma} \vz \sqrt h\right]) , ~~\vz \hbox{ being one realisation of a $\N_{0,1}^d$ r.v.}
\]
Bellman's dynamic programming principle tells us that the optimal control of the problem is the minimiser, $u^m$, in:
\eq{\label{DPP}&&
v^m(\vX):=
\min_{\vu\in\U} \Big\{
 \int_{m{h}}^{(m+1){h}}\E\left[|\bm\zeta_{{h}}(\vX,\vu,\vz)- \vX_d|^2-\bm\alpha\cdot\vu \right]\d\tau
\cr&& 
~~~~~~~~~~~~~~~+\beta[\vu]_t^{mh,(m+1)h}+\E[v^{m+1}(\bm\zeta_{{h}}(\vX,\vu,\vz)]
\Big\}
\cr&&
\approx\min_{\vu\in\U}  \left\{h\E\left[|\bm\zeta_{{h}}(\vX,\vu,\vz)- \vX_d|^2\right]-h\bm\alpha\cdot\vu + \beta\E\left[|\vu^{m+1}(\bm\zeta_{{h}}(\vX,\vu,\vz))- \vu|^2\right]
\right.\cr&&\left.
~~~~~ ~~~~ +\E[v^{m+1}(\bm\zeta_{{h}}(\vX,\vu,\vz)]\right\},\hbox{ with $v^M(\vX)=0.$}
}
Evidently,
\[
\E\left[|\bm\zeta_{{h}}(\vX,\vu,z)-\vX_d|^2\right]
=|\vX-\vX_d + {h}\vX\underline{\bm\Lambda}[\vr-\underline{\bm\kappa}\vX - \vu]|^2 + h|\underline{\bm\sigma}\vX|^2.
\]
For every component, to compute $\E[v^{m+1}(\bm\zeta_{{h}}(\vX,\vu,\vz))]$, we use a  quadrature formula with $Q$ points $\{z_q\}_1^Q$ and
weights $\{w_q\}_1^Q$ based on optimal quantization of the normal distribution $\N^d_{0,1}$ (see~\cite{PaPr2003} or~\cite[Chapter 5]{Pagesbook2018} and the website~\url{www.quantize.maths-fi.com} for download of grids) so that
\eq{\label{minquant}
	\E[v^{m+1}(\bm\zeta_{{h}}(\vX,\vu,\vz))]\approx  \sum_{q=1}^Q{w}_q v^{m+1}(\bm\zeta_{{h}}(\vX,\vu,\vz_q))
}
Finally at every time step and every $X_j=j L/J,$ $j=0, \dots, J$ (with $L>>1$), the result is minimised with respect to $\vu\in\U$ by a dichotomic search.
In this fashion $\{\vu^m(\vX_j)\}_{j=1}^J$ is obtained on a grid and a piecewise linear interpolation is constructed to prepare for the next time step $t^{m-1}$.

\subsection{Numerical results for a single species}\label{nrone}
For the numerical tests we have chosen: $r=2,\kappa=1.2, X_d=1, T=2$,
\[
\alpha=0.01,\quad \beta=0.1,\quad \sigma=0.1,\quad L=3, \quad J=40,\quad M=50, \quad K=100
\]
and $X_0=0.5+0.1\,j,\; j=0,\dots,9$.
Figures \ref{shjbex1} and \ref{resDP2} show the control surface $(X,t)\mapsto u(X,t)$ and the value function
$(X,t)\mapsto v(X,t)$.  Although the control seems to be either 0.5 or 1 everywhere, there is a small interval near $X=X_d=1$ in which it is not bang-bang. It is an important region, as seen
on the sample solutions which are shown on Figure \ref{s2-HJBSTO} and  \ref{s3-HJBSTO} where it is clear that the control is not always 0.5 or 1.
For these an initial condition $X_0$ is chosen, a noise $\d W$ is generated and the SDE are integrated by the Euler scheme with the approximately optimal control $u(X_t,t)$ obtained by the SDP method described above. Then, by definition of the cost function, $X_t$ should tend to be equal to $X_d$ as much as possible without too many jumps for $u$.  The trajectory is compared with one on which $u=u_M$, \textit{i.e.}, no quota.

Finally on Figure \ref{shjbex2} the cost function of the control problem is plotted for 100 samples and various values of $X_0$. Away from $X_d$, it increases, which makes sense because $X_0=1$ means that we start with the optimal value in terms of fish biomass. 
\begin{figure}[h!]
\begin{minipage} [b]{0.49\textwidth}
\centering
\tikzset{external/remake next}
\begin{tikzpicture}[scale=0.9,]%
\begin{axis}[legend style={at={(1,1)},anchor=north east}, compat=1.3,xlabel= {$X$},ylabel= {$t$}]
\addplot3[surf] table [ ] {fig/hjb_sto_results_u_trunc_M50.txt}; 
\addlegendentry{$u(X,t)$}
\end{axis}
\end{tikzpicture}
\caption{\label{shjbex1} Solution $u(X,t)$ from SDP.}
\end{minipage}
\hskip 0.5cm %
\begin{minipage} [b]{0.49\textwidth}
\begin{center}
\tikzset{external/remake next}
\begin{tikzpicture}[scale=0.8]
\begin{axis}[legend style={at={(1,1)},anchor=north east}, compat=1.3,xlabel= {$X_0$},ylabel= {average cost}, 
ytick={-0.01,-0.005,0,0.005,0.01,0.02,0.03,0.04,0.05}, grid style=dashed, ymajorgrids=true]
\addplot[thick,dotted,color=black,mark=none,mark size=1pt] table [x index=0, y index=1]{fig/hjb_sto_costs_M50_forlatex.txt};
\addlegendentry{$n=1$}
\addplot[thick,dashed,color=blue,mark=none,mark size=1pt] table [x index=0, y index=1]{fig/hjb_sto_costs_M100_forlatex.txt};
\addlegendentry{$n=2$}
\addplot[thick,solid,color=red,mark=square,mark size=1pt] table [x index=0, y index=1]{fig/hjb_sto_costs_M200_forlatex.txt};
\addlegendentry{$n=4$}
\end{axis}
\end{tikzpicture}
\caption{\label{shjbex2} Values of the average cost versus $X_0$ for $100$ realizations of $\d W$ and with the SDP, for the 3 time meshes $50\,n$, showing that the results are independent of $n$.}
\end{center}
\end{minipage}
\end{figure}
\begin{figure}[h!]
\begin{minipage} [b]{0.49\textwidth}
\centering
\tikzset{external/remake next}
\begin{tikzpicture}[scale=0.8]
\begin{axis}[legend style={at={(1,1)},anchor=north east}, compat=1.3, ymax=1.4,
  xlabel= {$time$},
  ylabel= {$u(X_t,t),~~~~~~X_t$}]
\addplot[thick,solid,color=red,mark=none,mark size=1pt] table [x index=0, y index=1]{fig/hjb_sto_results_traj_2.txt};
\addlegendentry{ $X_t$ with quota}
\addplot[thin,dashed,color=black,mark=none,mark size=1pt] table [x index=0, y index=2]{fig/hjb_sto_results_traj_2.txt};
\addlegendentry{ quota $u_t$}
\addplot[thin,solid,color=blue,mark=none,mark size=1pt] table [x index=0, y index=3]{fig/hjb_sto_results_traj_2.txt};
\addlegendentry{ $X_t$ without quota}
\end{axis}
\end{tikzpicture}
\caption{\label{s2-HJBSTO} Simulation of the fishing model with a quota function computed by SDP and $X_0=0.7$.} 
\end{minipage}
\hskip 0.5cm
\begin{minipage} [b]{0.49\textwidth}
\centering
\tikzset{external/remake next}
\begin{tikzpicture}[scale=0.8]
\begin{axis}[legend style={at={(1,1)},anchor=north east}, compat=1.3, ymax=1.4,
  xlabel= {$time$},
  ylabel= {$u(X_t,t),~~~~~~X_t$}]
\addplot[thick,solid,color=red,mark=none,mark size=1pt] table [x index=0, y index=1]{fig/hjb_sto_results_traj_8.txt};
\addlegendentry{ $X_t$ with quota}
\addplot[thin,dashed,color=black,mark=none,mark size=1pt] table [x index=0, y index=2]{fig/hjb_sto_results_traj_8.txt};
\addlegendentry{ quota $u_t$}
\addplot[thin,solid,color=blue,mark=none,mark size=1pt] table [x index=0, y index=3]{fig/hjb_sto_results_traj_8.txt};
\addlegendentry{ $X_t$ without quota}
\end{axis}
\end{tikzpicture}
\caption{Simulation of the fishing model with a quota function computed by SDP and $X_0=1.3$.}
\label{s3-HJBSTO} 
\end{minipage}
\end{figure}

\section{Hamilton-Jacobi-Bellman solutions (HJB)}
Let us return to (\ref{DPP}) and note that:
\begin{eqnarray}\label{calc}&
v^m(\vX)\approx \min_{\vu\in\U}  \Big\{ |\vX - \vX_d|^2 
-\bm\alpha\cdot\vu &
+ \frac\beta h\E\left[|\vu^{m+1}(\bm\zeta_{{h}}(\vX,\vu,\vz))- \vu(\vX)|^2\right]
\cr&&
+\E[v^{m+1}({\bm\zeta}_{{h}}(\vX,\vu,\vz)]\Big\},
\end{eqnarray}
with $v^M=0$.  According to Remark \ref{remm1}
\eq{&&
\E\left[|\vu^{m+1}(\bm\zeta_{{h}}(\vX,\vu,\vz))- \vu|^2\right]
 \approx 
 \E\left[\left|\vu^{m+1}\left(\vX+ \vX\underline{\bm\Lambda}\left[h(\vr-\underline{\bm\kappa}\vX -\vu) +\sqrt h \underline{\bm\sigma}\vz\right]\right)-\vu\right|^2\right]
 \cr &&
\approx \left|\vu^{m+1}\left(\vX+ \vX\underline{\bm\Lambda}\left[h(\vr-\underline{\bm\kappa}\vX)\right]\right) -\vu\right|^2 + \frac h2|\underline{\bm\sigma}\vX\vu'^{m+1}|^2.
}
Also, by It\^o formula,
\eq{&&
\E[v^{m+1}(\zeta_{{h}}(\vX,\vu^m,\vz)] = \E[v^{m+1}(\vX+h\vX\underline{\bm\Lambda}[\vr-\kappa \vX - \vu^{m}]+\vX\underline{\bm\sigma}\sqrt{h}\N^d_{0,1})]
\cr&&
=
v^{m+1}(\vX+h\vX\underline{\bm\Lambda}[\vr-\underline{\bm\kappa} \vX - \vu^{m}])+\frac{h}2 |\underline{\bm\sigma}\vX|^2{v''}^{m+1}(\vX).
}
Consequently, with $v''^{m+1}$ replaced by $v''^m$ for numerical convenience,  (\ref{DPP})  becomes
\eq{\label{hjb1}&&
v^m(\vX)\approx \min_{\vu^m\in\U}\Big\{v^{m+1}(\vX+h\vX\underline{\bm\Lambda}[\vr-\underline{\bm\kappa} \vX - \vu^{m}(\vX))] + \frac{h}2 |\underline{\bm\sigma}\vX|^2{v''}^{m+1}(\vX)
+ h|\vX - \vX_d|^2
\cr&&
-h\bm\alpha\cdot\vu^m(\vX) + \beta|\vu^{m+1}(\vX+h\vX\underline{\bm\Lambda}[\vr-\underline{\bm\kappa} \vX-\vu^{m})]- \vu^m(\vX)|^2
+ h\frac{\beta}{2}|\underline{\bm\sigma}\vX\vu'^{m+1}|^2\Big\}.
}
The last line is an approximation motivated by the search for a stable implicit numerical scheme for $v^m$:
\eq{\label{hjbd}&&
 \frac1h v^m(\vX) + \frac{|\underline{\bm\sigma}\vX|^2}2{v''}^{m}(\vX)
 =\frac1h v^{m+1}(\vX+h\vX\underline{\bm\Lambda}[\vr-\underline{\bm\kappa} \vX - \vu^{m}])  
+ |\vX - \vX_d|^2-\bm\alpha\cdot \vu^m(\vX) 
\cr&&
+ \frac\beta h|\vu^{m+1}(\vX+h\vX\underline{\bm\Lambda}[\vr-\underline{\bm\kappa} \vX-\vu^{m}])- \vu^m(\vX)|^2+\frac{\beta}2|\underline{\bm\sigma\vX}^2|\vu'^{m+1}|^2.
}
Furthermore, $\vu^m$ is given by minimising the right hand side of (\ref{hjbd}): 
\eq{\label{fp}
\vu^m =  \vu^{m+1}(\bm\zeta_h(\vX,\vu^m,0) +\frac{h}{2\beta}(\bm\alpha+\underline{\bm\Lambda}[\vX] \n v^{m+1}(\eta_h(\vX,\vu^m,0)) + o(h,\frac h\beta),
}
capped by the bounds $u_m$ and $u_M$.
Boundary conditions are not needed at $\vX=0$ because the PDE is self starting but for large $\vX$, cancelation of the quadratic terms requires ${v''}^m(\vX)\sim -2/||\underline{\bm\sigma}||^2$.
\begin{remark}
Notice that (\ref{fp}) makes an important use of the comment in Remark \ref{remm1}.
\end{remark}
\subsection{Numerical results for a single species by HJB}

The finite element method of degree~1 on intervals was used to approximate spatially~\eqref{hjbd}.  The linear systems were solved using the FreeFem++ software~\cite{freefem}. The range of $\vX$ was approximated by the segment $(0,3)$ discretized into $160$ intervals; $200$ time steps were used.

With the same parameters as in Section \ref{nrone} the results 
  are shown on Figures  \ref{hjbu}, \ref{resDP1}, \ref{simul2}, \ref{simul3}  and \ref{resDP3}.  These should be compared, respectively, with Figures \ref{shjbex1}, \ref{resDP2}, \ref{s2-HJBSTO},  \ref{s3-HJBSTO} and \ref{shjbex2}.

Except for the singularity at final time, which is not handled in the same way, the results are similar. HJB seems to handle better the penalty term with $\beta$ as $t\mapsto u(t)$ is not bang-bang.

\begin{figure}[h!]
\begin{minipage} [b]{0.49\textwidth}
\centering
\tikzset{external/remake next}
\begin{tikzpicture}[scale=0.9]
\begin{axis}[legend style={at={(1,1)},anchor=north east}, compat=1.3,xlabel= {$X$},ylabel= {$t$}]
 \addplot3[surf] table [ ] {fig/hjb_sto_results_V_M50.txt};
\addlegendentry{ $V(X,t)$}
\end{axis}
\end{tikzpicture}
\caption{SDP solution: $V(X,t)$.}
\label{resDP2} 
\end{minipage}
\begin{minipage} [b]{0.49\textwidth}
\centering
\tikzset{external/remake next}
\begin{tikzpicture}[scale=0.9]
\begin{axis}[legend style={at={(1,1)},anchor=north east}, compat=1.3,xlabel= {$X$},ylabel= {$t$}]
 \addplot3[surf] table [ ] {fig/bellman2.txt};
\addlegendentry{ $V(X,t)$}
\end{axis}
\end{tikzpicture}
\caption{HJB solution: value function $V(X,t)$.}
\label{resDP1} 
\end{minipage}\hskip 0.5cm
\end{figure}

\begin{figure}[h!]
\begin{minipage} [b]{0.49\textwidth}
\centering
\tikzset{external/remake next}
\begin{tikzpicture}[scale=0.9,]%
\begin{axis}[legend style={at={(1,1)},anchor=north east}, compat=1.3,xlabel= {$X$},ylabel= {$t$}]
 \addplot3[surf] table [ ] {fig/bellman1.txt}; 
\addlegendentry{ $u(X,t)$}
\end{axis}
\end{tikzpicture}
\caption{\label{hjbu} Solution $u(X,t)$ from HJB.}
\end{minipage}\hskip 0.5cm
\begin{minipage} [b]{0.45\textwidth}
\begin{center}
\tikzset{external/remake next}
\begin{tikzpicture}[scale=0.8]
\begin{axis}[legend style={at={(1,1)},anchor=north east}, compat=1.3,xlabel= {$X_0$},ylabel= {average cost}, 
ytick={-0.01,-0.005,0,0.005,0.01,0.02,0.03,0.04,0.05}, grid style=dashed, ymajorgrids=true]
\addplot[thick,dotted,color=black,mark=none,mark size=1pt] table [x index=0, y index=1]{fig/errcurv.txt};
\addlegendentry{$n=1$}
\addplot[thick,dashed,color=blue,mark=none,mark size=1pt] table [x index=0, y index=2]{fig/errcurv.txt};
\addlegendentry{$n=2$}
\addplot[thick,solid,color=red,mark=square,mark size=1pt] table [x index=0, y index=3]{fig/errcurv.txt};
\addlegendentry{$n=4$}
\end{axis}
\end{tikzpicture}
\caption{\label{resDP3}Values of the cost functions for the 3 space$\times$times meshes $40n\times50n$, to solve the control problem by a HJB algorithm.}
\end{center}
\end{minipage} 
\end{figure}
\begin{figure}[h!]
\begin{minipage} [b]{0.49\textwidth}
\centering
\tikzset{external/remake next}
\begin{tikzpicture}[scale=0.8]
\begin{axis}[legend style={at={(1,1)},anchor=north east}, compat=1.3, ymax=1.4,
  xlabel= {$time$},
  ylabel= {$u(X_t,t),~~~~~~X_t$}]
\addplot[thick,solid,color=red,mark=none,mark size=1pt] table [x index=0, y index=1]{fig/bellmanx2.txt};
\addlegendentry{ $X_t$ with quota}
\addplot[thin,dashed,color=black,mark=none,mark size=1pt] table [x index=0, y index=2]{fig/bellmanx1.txt};
\addlegendentry{ quota $u_t$}
\addplot[thin,solid,color=blue,mark=none,mark size=1pt] table [x index=0, y index=3]{fig/bellmanx2.txt};
\addlegendentry{ $X_t$ without quota}
\end{axis}
\end{tikzpicture}
\caption{Simulation of the fishing model with a quota function computed by HJB and $X_0=0.7$.}
\label{simul2} 
\end{minipage}
\hskip 0.5cm
\begin{minipage} [b]{0.49\textwidth}
\centering
\tikzset{external/remake next}
\begin{tikzpicture}[scale=0.8]
\begin{axis}[legend style={at={(1,1)},anchor=north east}, compat=1.3, ymax=1.4,
  xlabel= {$time$},
  ylabel= {$u(X_t,t),~~~~~~X_t$}]
\addplot[thick,solid,color=red,mark=none,mark size=1pt] table [x index=0, y index=1]{fig/bellmanx8.txt};
\addlegendentry{ $X_t$ with quota}
\addplot[thin,dashed,color=black,mark=none,mark size=1pt] table [x index=0, y index=2]{fig/bellmanx8.txt};
\addlegendentry{ quota $u_t$}
\addplot[thin,solid,color=blue,mark=none,mark size=1pt] table [x index=0, y index=3]{fig/bellmanx8.txt};
\addlegendentry{ $X_t$ without quota}
\end{axis}
\end{tikzpicture}
\caption{Simulation of the fishing model with a quota function computed by HJB and $X_0=1.3$.}
\label{simul3} 
\end{minipage}
\end{figure}

\section{Optimisation with Markovian feedback neural network controls}

In this section, we look for Markovian feedback controls, \textit{i.e.}, functions $u:[0,T] \times \R^+ \to \R$, adapted to $X_t$, satisfying $u_m\leq u(x,t)\leq u_M$ for every $(x,t)$. We restrict our attention to such functions which are encoded by neural networks. This strategy has been used previously in the optimal control literature, see \textit{e.g.}~\cite{GobetMunos05},\cite{HanE16},\cite{CarmonaLauriere19b},\cite{PHA}.

Given an {\it activation} function $\psi:\R\to\R$ (\textit{e.g.} a sigmoid or \texttt{ReLU}), let us introduce some helpful notation:
\eq{
	\mathbf{L}_{d_1, d_2} = \left\{\phi: z\in\R^{d_1}\mapsto\{ \psi(\beta_i + w_i\cdot z) \}_{i=1}^{d_2}\in\R^{d_2}\,\Big|\,   \beta \in \RR^{d_2},  w \in \RR^{d_2 \times d_1} \right\} 
}
is the set of layer functions with input dimension $d_1$ and output dimension $d_2$. 

\begin{figure}
\centering
\includegraphics[width=10cm]{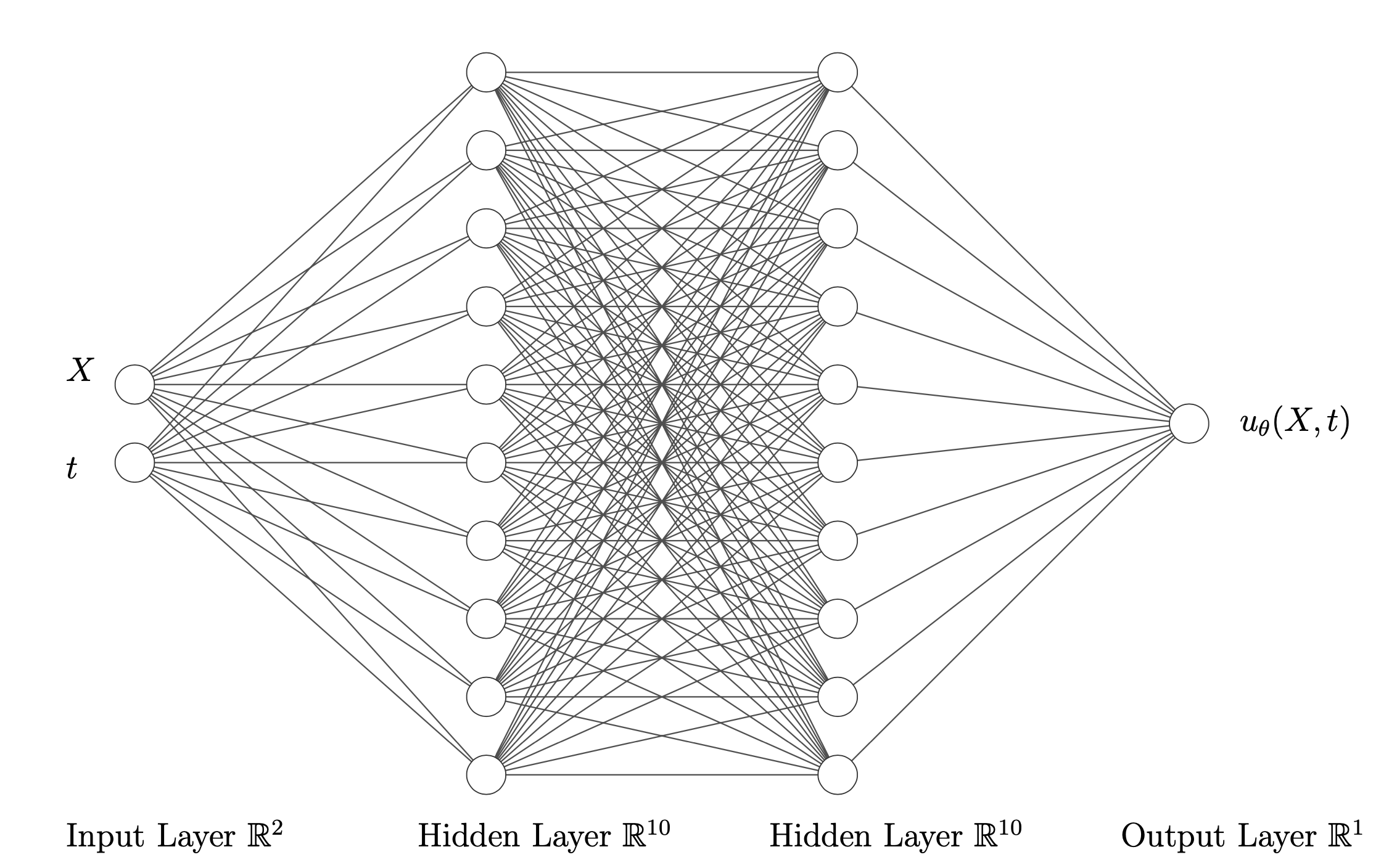}
\caption{\label{sketch}Sketch of a two layers neural network with 2 inputs, $X,t$ and one output $u$.  Each layer here has 10 neurons while for this article we used 100 neurons in each layers.}
\end{figure}

Let $\vd=\{d_0, \dots, d_{\ell+1}\}$; define  the set of feedforward fully connected neural networks with $\ell$ hidden layers and one output layer,
\eq{\label{nnd}
	\mathbf{N}_{\vd} 
	&= 
	\left\{ u: \RR^{d_0} \to \RR^{d_{\ell+1}},\, u = \phi^{(\ell)} \circ \phi^{(\ell-1)} \circ \dots \circ \phi^{(0)} \,\Big|\,  \phi^{(i)} \in \mathbf{L}_{d_{i}, d_{i+1}},\, i=0,\dots,\ell
 \right\} 
}
where $\circ$ denotes the composition of functions.
Denote by $\Theta$ the set of parameters:
 \[
 \Theta =\left\{\,\theta:=(\beta^{(0)}, w^{(0)},\beta^{(1)}, w^{(1)},\cdots,\beta^{(\ell)}, w^{(\ell)})\Big| \quad \beta^{(i)} \in \RR^{d_{i}},  w^{(i)} \in \RR^{d_{i} \times d_{1+1}} \right\} .
\]
 For each $\theta\in\Theta$, the corresponding network function will be denoted by  $\vu_\theta$. Hence (\ref{nnd}) may be rewritten as:
\[
\mathbf{N}_{\vd}  = \{\vu_\theta \,:\, \theta \in \Theta\}
\]
For the fishing control problem  $d_0 = d+1$ ($d$ species for the space variable, plus the time variable) and $d_{\ell+1} = d$ (the dimension of the control variable $\vu$). Furthermore, to ensure the constraint $\vu_m\leq \vu(x,t)\leq \vu_M$, it is sufficient to take the last activation function which satisfies this constraint. In the present implementation we have used $\vu_m + (\vu_M-\vu_m)\sigma(x)$ where $\sigma$ is the sigmoid function.  Figure \ref{sketch} shows a fully connected 2 layers neural network with 2 inputs, one output and 10 neurons in each layer.

The control problem is approximated by the minimisation over $\theta\in\Theta$ of:
\eq{\label{NN}
\mathbb{J}(\theta)=\frac hK\sum_{k=1}^K\sum_{m=1}^M \|\vX^m_k-\vX_d\|^2-\bm\alpha\cdot \vu_\theta(\vX_k^m,t^m)
+\frac\beta{h}|\vu_\theta(\vX_k^m,t^m) - \vu_\theta(\vX_k^{m-1},t^{m-1})|^2
}
where $\vX_k^m$ is the result of one step of the Euler scheme \eqref{euler1} with $\vu_\theta(\vX_k^{m-1},t^{m-1})$ and a  realization of the noise. 

Noting that (\ref{NN}) is a sum of terms, we can run a Stochastic Gradient Descent or one of its variants like \texttt{ADAM}. At each iteration, we sample a mini-batch of initial positions and realizations of  the noise which allow us to compute random realizations of trajectories and compute a partial sum of (\ref{NN}). It is then used to compute a gradient and back-propagate it to adjust $\theta$.  

The stochastic gradient algorithm \texttt{ADAM} is used to update $\theta$ towards a local minimum of $\mathbb{J}(\theta)$ needs the gradient of $\mathbb{J}(\theta)$.  A mini-batch is a random subset $B\subset\{1,\dots,K\}$; the gradients with respect to $\theta$ are computed as follows. Let  $\mathbb{J}_B(\theta) :=\frac hK\sum_{k\in B} \textbf{j}_k(\theta)$ where  
\[
\textbf{j}_k(\theta)=\sum_{m=1}^M \|\vX^m_k-\vX_d\|^2-\bm\alpha\cdot \vu_\theta(\vX_k^m,t^m)
+\frac\beta{h}|\vu_\theta(\vX_k^m,t^m) - \vu_\theta(\vX_k^{m-1},t^{m-1})|^2.
\]
The gradient of $\textbf{j}_k(\theta)$ with respect to $\theta=\{\theta_n\}_n$ is made of: %
\eq{
\frac{\d\textbf{j}_k}{\d \theta_n} =\sum_m\frac{\d\textbf{j}_k}{\d \vu^m}\frac{\d \vu^m}{\d\theta_n} \hbox { with }\frac{\d\textbf{j}_k}{\d \vu^m}= P^m_k\vX^m_k -\bm\alpha - 2\frac\beta{h}(\vu^{m+1}-2\vu^m+\vu^{m-1}),
}
where $P^m_k$ is defined by $P^{M-2}_k=0$ and
\eq{
P^{m-1}_k= P^m_k \left(1+{h}(r-2 \kappa \vX^m_k -\vu^m) + \sigma\d W^m_k,\right)- 2(\vX^m_k-\vX_d) ,~~m=M-2,\dots,1.
}
 
Softwares for neural networks such as \texttt{tensorflow} provide automatic differentiation tools (back propagation) to compute $\frac{\d \vu^m}{\d\theta_n}$.

Recall that \texttt{ADAM} algorithm consists in choosing $B$ randomly and then
\begin{enumerate}
\item Let $\alpha=0.001$, $m_1=0.9$, $m_2=0.999$, $\varepsilon=10^{-8}$.

 Loop on i:
 
\item Set ${\bf g}^i=\{\frac{\d\mathbb{J}_B}{\d \theta^i_n}\}_{n=1}^N$
\item Set ${\bf p}^i=m_1 {\b p}^{i-1} + (1-m_1){\bf g}^i$
\item Set $q^i= m_2 q^{i-1} + (1-m_2)|{\bf g^i}|^2$
\item Set $m_j^i=m_j\cdot m_j^{i-1},\, j=1,2$,\quad $\ds \hat{\bf p}^i=\frac{{\bf p}^i}{1-m_1^i}$, $\ds \hat q^i=\frac{q^i}{1-m_2^i}$ 
\item Update  $\ds \theta^i_n= \theta^{i-1}_n -\frac{\alpha\cdot \hat p^i_n}{\sqrt{\hat q^i} + \varepsilon}$, $ n=1,\dots,N$, the number of parameters.
\end{enumerate}

\subsection{ Numerical results using a neural network for a single species}

Notice that we use the tools of AI but there is no learning phase from known solutions, rather the learning phase is simply to find the best parameters to achieve the minimum of a given cost function.

To check the results we implemented also a single layer neural network with 15000 neurons directly in C++ using the automatic differentiation in reverse mode of the library \texttt{adept} \cite{adept} and a conjugate gradient algorithm instead of ADAM.  Convergence is fast (20 to 40 iterations) but this simple method does not work for multiple layered neural network because of memory limitation.  

For the single species numerical test, the results are the same whether obtained by the one layer or the two layers networks.

Results are shown on Figures \ref{mathieu3}, \ref{simul2-HJBSTO}, \ref{simul3-HJBSTO} and \ref{mathieu4}
and should be compared, respectively, with Figures  \ref{hjbu}, \ref{simul2}, \ref{simul3}  and \ref{resDP3} and/or, respectively, with Figures \ref{shjbex1}, \ref{s2-HJBSTO},  \ref{s3-HJBSTO} and \ref{shjbex2}.

Note that the results are closer to those obtained with HJB than those obtained with SDP.

Our general impression is that all 3 methods are equal, with a slight reserve for SDP for which the control oscillates much more than with the other two methods, implying that it does not handle as well the $\beta$-penalisation term.
\begin{figure}[h!]
\begin{minipage} [b]{0.49\textwidth}
\centering
\tikzset{external/remake next}
\begin{tikzpicture}[scale=0.9]
\begin{axis}[legend style={at={(1,1)},anchor=north east}, compat=1.3,xlabel= {$X$},ylabel= {$t$}]
 \addplot3[surf] table [ ] {fig/nn_control_Nt50.txt};
\addlegendentry{ $u_\theta(X,t)$}
\end{axis}
\end{tikzpicture}
\caption{\label{mathieu3} Dynamic feedback control,$u_\theta(X,t)$ computed by the Markovian Neural Network.} 
\end{minipage}
										\hskip 1.0cm
\begin{minipage} [b]{0.45\textwidth}
\begin{center}
\tikzset{external/remake next}
\begin{tikzpicture}[scale=0.8]
\begin{axis}[legend style={at={(1,1)},anchor=north east}, compat=1.3,xlabel= {$X_0$},ylabel= {average cost}, 
ytick={-0.01,-0.005,0,0.005,0.01,0.02,0.03,0.04,0.05}, grid style=dashed, ymajorgrids=true]
\addplot[thick,dotted,color=black,mark=none,mark size=1pt] table [x index=0, y index=1]{fig/nn_costs_Nt50.txt};
\addlegendentry{$n=1$}
\addplot[thick,dashed,color=blue,mark=none,mark size=1pt] table [x index=0, y index=1]{fig/nn_costs_Nt100.txt};
\addlegendentry{$n=2$}
\addplot[thick,solid,color=red,mark=square,mark size=1pt] table [x index=0, y index=1]{fig/nn_costs_Nt200.txt};
\addlegendentry{$n=4$}
\end{axis}
\end{tikzpicture}
\caption{\label{mathieu4} Single species: Values of the  cost versus $X_0$ for $100$ realizations with the Markovian neural network, for the 3 time meshes $50\times n$.}
\end{center}
\end{minipage} 
\end{figure}

\begin{figure}[h!]
\begin{minipage} [b]{0.49\textwidth}
\centering
\tikzset{external/remake next}
\begin{tikzpicture}[scale=0.8]
\begin{axis}[legend style={at={(1,1)},anchor=north east}, compat=1.3, ymax=1.4,
  xlabel= {$time$},
  ylabel= {$u(X_t,t),~~~~~~X_t$}]
\addplot[thick,solid,color=blue,mark=none,mark size=1pt] table [x index=0, y index=1]{fig/nn_traj2_Nt50.txt};
\addlegendentry{ $X_t$ with quota}
\addplot[thin,dashed,color=black,mark=none,mark size=1pt] table [x index=0, y index=2]{fig/nn_traj2_Nt50.txt};
\addlegendentry{ quota $u_t$}
\addplot[thin,solid,color=red,mark=none,mark size=1pt] table [x index=0, y index=3]{fig/nn_traj2_Nt50.txt};
\addlegendentry{ $X_t$ without quota}
\end{axis}
\end{tikzpicture}
\caption{\label{simul2-HJBSTO} Simulation for a single fish species computed by the Markovian Neural Network and $X_0=0.7$. Performance with and without quota $u_t$.} 
\end{minipage}
\hskip 0.5cm
\begin{minipage} [b]{0.49\textwidth}
\centering
\tikzset{external/remake next}
\begin{tikzpicture}[scale=0.8]
\begin{axis}[legend style={at={(1,1)},anchor=north east}, compat=1.3, ymax=1.4,
  xlabel= {$time$},
  ylabel= {$u(X_t,t),~~~~~~X_t$}]
\addplot[thick,solid,color=blue,mark=none,mark size=1pt] table [x index=0, y index=1]{fig/nn_traj8_Nt50.txt};
\addlegendentry{ $X_t$ with quota}
\addplot[thin,dashed,color=black,mark=none,mark size=1pt] table [x index=0, y index=2]{fig/nn_traj8_Nt50.txt};
\addlegendentry{ quota $u_t$}
\addplot[thin,solid,color=red,mark=none,mark size=1pt] table [x index=0, y index=3]{fig/nn_traj8_Nt50.txt};
\addlegendentry{ $X_t$ without quota}
\end{axis}
\end{tikzpicture}
\caption{\label{simul3-HJBSTO} Simulation for a single fish species computed by the Markovian Neural Network and $X_0=1.3$. Performance with and without quota $u_t$.}
\end{minipage}
\end{figure}

\section{Numerical results: 3 species}

We now turn our attention to an example with 3 species. Consider the following interaction matrix between species:
\[
\kappa = \left(\begin{matrix}
1.2&-0.1&0\cr 0.2&1.2&0 \cr 0.1&0.1&1\cr
\end{matrix}\right)
\]
All other parameters are as in Section \ref{nrone} and $X_d=(1,1,1)^T$.

With SDP (Stochastic Dynamic Programming) there is a difficulty: the minimisation in (\ref{minquant}) is now 3 dimensional and cannot be obtained by dichotomy.  So the method was not tested.  We have compared HJB and Neural Network optimisation with 1 and then with 2 layers.  The results differ; according to the last subsection none of the methods found the true solution.

\subsection{Numerical results for 3 species with HJB}
The computation is done with $80$ time steps. Here too a finite element discretisation was used with $P^1$ elements on tetraedra. The mesh for the cube $(0,3)^3$ is obtained from an automatic mesh generator from a $30\times 30\times 30)$ surface mesh resulting into 29791 vertices and 16200 tetraedra. 
It took 7 min on an M1 Apple laptop. 
The vector valued optimal control $\vX,t\in\R^3\times(0,T)\mapsto\vu\in(0.5,1)^3$ computed with HJB is shown on Figures \ref{bellman11} to \ref{bellman33}. 
In these, two of the 3  coordinates of $\vX$, are fixed at their $X_d$ values.

Then the fishing model is integrated with the optimal $\vu$ and a random realization of $\d\vW$ for 2 values of $\vX_0$: $\vX_0=(0.7,0.7,0.7)^T$ or $(1.3,1.3,1.3)^T$.  The results are compared with a similar simulation without quota, \textit{i.e.}, $\vu\equiv 1$.  See Figures
\ref{simul3hjb1}, \ref{simul3hjb2}, \ref{simul3hjb3}.
The quality of the optimisation is seen visually when $X_t$ is closest to $X_d$ and numerically from the lowest value of the cost function, plotted versus $\vX_0$ on Figure \ref{hjb3cost}.
\begin{figure}[h!]
\tikzset{external/remake next}
\begin{minipage} [b]{0.32\textwidth}
\centering
\begin{tikzpicture}[scale=0.6]
\begin{axis}[legend style={at={(1,1)},anchor=north east}, compat=1.3,xlabel= {$X$},ylabel= {$t$}]
 \addplot3[surf] table [ ] {results/hjb3b/bellman11.txt};
\addlegendentry{ $u_1(X_1,t)$}
\end{axis}
\end{tikzpicture}
\caption{\label{bellman11} HJB:  $X_1,t\mapsto u_1$.}
\end{minipage}\hskip 0.5cm
\begin{minipage} [b]{0.32\textwidth}
\centering
\begin{tikzpicture}[scale=0.6]
\begin{axis}[legend style={at={(1,1)},anchor=north east}, compat=1.3,xlabel= {$X$},ylabel= {$t$}]
 \addplot3[surf] table [ ] {results/hjb3b/bellman12.txt};
\addlegendentry{ $u_1(X_2,t)$}
\end{axis}
\end{tikzpicture}
\caption{\label{bellman12} HJB: $X_2,t\mapsto u_1$.} 
\end{minipage}
\begin{minipage} [b]{0.32\textwidth}
\centering
\begin{tikzpicture}[scale=0.6]
\begin{axis}[legend style={at={(1,1)},anchor=north east}, compat=1.3,xlabel= {$X$},ylabel= {$t$}]
 \addplot3[surf] table [ ] {results/hjb3b/bellman13.txt};
\addlegendentry{ $u_1(X_3,t)$}
\end{axis}
\end{tikzpicture}
\caption{\label{bellman13}HJB: $X_3,t\mapsto u_1$.} 
\end{minipage}
\end{figure}

\begin{figure}[h!]
\tikzset{external/remake next}
\begin{minipage} [b]{0.32\textwidth}
\centering
\begin{tikzpicture}[scale=0.6]
\begin{axis}[legend style={at={(1,1)},anchor=north east}, compat=1.3,xlabel= {$X$},ylabel= {$t$}]
 \addplot3[surf] table [ ] {results/hjb3b/bellman21.txt};
\addlegendentry{ $u_2(X_1,t)$}
\end{axis}
\end{tikzpicture}
\caption{HJB: $X_1,t\mapsto u_2$.}
\label{bellman21} 
\end{minipage}\hskip 0.5cm
\begin{minipage} [b]{0.32\textwidth}
\centering
\begin{tikzpicture}[scale=0.6]
\begin{axis}[legend style={at={(1,1)},anchor=north east}, compat=1.3,xlabel= {$X$},ylabel= {$t$}]
 \addplot3[surf] table [ ] {results/hjb3b/bellman22.txt};
\addlegendentry{ $u_2(X_2,t)$}
\end{axis}
\end{tikzpicture}
\caption{HJB: $X_2,t\mapsto u_2$.}
\label{bellman22} 
\end{minipage}
\begin{minipage} [b]{0.32\textwidth}
\centering
\begin{tikzpicture}[scale=0.6]
\begin{axis}[legend style={at={(1,1)},anchor=north east}, compat=1.3,xlabel= {$X$},ylabel= {$t$}]
 \addplot3[surf] table [ ] {results/hjb3b/bellman23.txt};
\addlegendentry{ $u_2(X_3,t)$}
\end{axis}
\end{tikzpicture}
\caption{HJB: $X_3,t\mapsto u_2$.}
\label{bellman23} 
\end{minipage}
\end{figure}

\begin{figure}[h!]
\tikzset{external/remake next}
\begin{minipage} [b]{0.32\textwidth}
\centering
\begin{tikzpicture}[scale=0.6]
\begin{axis}[legend style={at={(1,1)},anchor=north east}, compat=1.3,xlabel= {$X$},ylabel= {$t$}]
 \addplot3[surf] table [ ] {results/hjb3b/bellman31.txt};
\addlegendentry{ $u_3(X_1,t)$}
\end{axis}
\end{tikzpicture}
\caption{HJB: $X_1,t\mapsto u_3$.}
\label{bellman31} 
\end{minipage}\hskip 0.5cm
\begin{minipage} [b]{0.32\textwidth}
\centering
\begin{tikzpicture}[scale=0.6]
\begin{axis}[legend style={at={(1,1)},anchor=north east}, compat=1.3,xlabel= {$X$},ylabel= {$t$}]
 \addplot3[surf] table [ ] {results/hjb3b/bellman32.txt};
\addlegendentry{HJB: $u_3(X_2,t)$}
\end{axis}
\end{tikzpicture}
\caption{HJB:$X_2,t\mapsto u_3$.}
\label{bellman32} 
\end{minipage}
\begin{minipage} [b]{0.32\textwidth}
\centering
\begin{tikzpicture}[scale=0.6]
\begin{axis}[legend style={at={(1,1)},anchor=north east}, compat=1.3,xlabel= {$X$},ylabel= {$t$}]
 \addplot3[surf] table [ ] {results/hjb3b/bellman33.txt};
\addlegendentry{ $u_3(X_3,t)$}
\end{axis}
\end{tikzpicture}
\caption{HJB:$X_3,t\mapsto u_3$.}
\label{bellman33} 
\end{minipage}
\end{figure}
\begin{figure}[h!]
\tikzset{external/remake next}
\begin{minipage} [b]{0.32\textwidth}
\centering
\begin{tikzpicture}[scale=0.6]
\begin{axis}[legend style={at={(1,1)},anchor=north east}, compat=1.3,
  xlabel= {$time$},
  ylabel= {$u(X_t,t),~~~~~~X_t$}]
\addplot[thick,solid,color=red,mark=none,mark size=1pt] table [x index=0, y index=1]{results/hjb3b/bellmanx2.txt};
\addlegendentry{ $X_1$}
\addplot[thin,solid,color=black,mark=none,mark size=1pt] table [x index=0, y index=2]{results/hjb3b/bellmanx2.txt};
\addlegendentry{  $u_1$}
\addplot[thick,solid,color=blue,mark=none,mark size=1pt] table [x index=0, y index=3]{results/hjb3b/bellmanx2.txt};
\addlegendentry{ $Y_1$}

\addplot[thick,dashed,color=red,mark=none,mark size=1pt] table [x index=0, y index=1]{results/hjb3b/bellmanx7.txt};
\addlegendentry{ $X_1$}
\addplot[thin,dashed,color=black,mark=none,mark size=1pt] table [x index=0, y index=2]{results/hjb3b/bellmanx7.txt};
\addlegendentry{ $u_1$}

\addplot[thick,dashed,color=blue,mark=none,mark size=1pt] table [x index=0, y index=3]{results/hjb3b/bellmanx7.txt};
\addlegendentry{ $Y_1$}

\end{axis}
\end{tikzpicture}
\caption{HJB:Optimal biomass and quota function when $X_1(0)=0.7$ and $X_1(0)=1.3$.}
\label{simul3hjb1} 
\end{minipage}
\begin{minipage} [b]{0.32\textwidth}
\tikzset{external/remake next}
\centering
\begin{tikzpicture}[scale=0.6]
\begin{axis}[legend style={at={(1,1)},anchor=north east}, compat=1.3,
  xlabel= {$time$},
  ylabel= {$u(X_t,t),~~~~~~X_t$}]
\addplot[thick,solid,color=red,mark=none,mark size=1pt] table [x index=0, y index=1]{results/hjb3b/bellmany2.txt};
\addlegendentry{ $X_2$}
\addplot[thin,solid,color=black,mark=none,mark size=1pt] table [x index=0, y index=2]{results/hjb3b/bellmany2.txt};
\addlegendentry{  $u_2$}
\addplot[thick,solid,color=blue,mark=none,mark size=1pt] table [x index=0, y index=3]{results/hjb3b/bellmany2.txt};
\addlegendentry{ $Y_2$}
\addplot[thick,dashed,color=red,mark=none,mark size=1pt] table [x index=0, y index=1]{results/hjb3b/bellmany7.txt};
\addlegendentry{ $X_2$}
\addplot[thin,dashed,color=black,mark=none,mark size=1pt] table [x index=0, y index=2]{results/hjb3b/bellmany7.txt};
\addlegendentry{ $u_2$}
\addplot[thick,dashed,color=blue,mark=none,mark size=1pt] table [x index=0, y index=3]{results/hjb3b/bellmany7.txt};
\addlegendentry{ $Y_2$}
\end{axis}
\end{tikzpicture}
\caption{HJB:Optimal biomass and quota function when $X_2(0)=0.7$ and $X_2(0)=1.3$.}
\label{simul3hjb2} 
\end{minipage}
\begin{minipage} [b]{0.32\textwidth}
\centering
\begin{tikzpicture}[scale=0.6]
\begin{axis}[legend style={at={(1,1)},anchor=north east}, compat=1.3,
  xlabel= {$time$},
  ylabel= {$u(X_t,t),~~~~~~X_t$}]
\addplot[thick,solid,color=red,mark=none,mark size=1pt] table [x index=0, y index=1]{results/hjb3b/bellmanz2.txt};
\addlegendentry{ $X_3$}
\addplot[thin,solid,color=black,mark=none,mark size=1pt] table [x index=0, y index=2]{results/hjb3b/bellmanz2.txt};
\addlegendentry{  $u_3$}
\addplot[thick,solid,color=blue,mark=none,mark size=1pt] table [x index=0, y index=3]{results/hjb3b/bellmanz2.txt};
\addlegendentry{ $Y_3$}
\addplot[thick,dashed,color=red,mark=none,mark size=1pt] table [x index=0, y index=1]{results/hjb3b/bellmanz7.txt};
\addlegendentry{ $X_3$}
\addplot[thin,dashed,color=black,mark=none,mark size=1pt] table [x index=0, y index=2]{results/hjb3b/bellmanz7.txt};
\addlegendentry{ $u_3$}
\addplot[thick,dashed,color=blue,mark=none,mark size=1pt] table [x index=0, y index=3]{results/hjb3b/bellmanz7.txt};
\addlegendentry{ $Y_3$}
\end{axis}
\end{tikzpicture}
\caption{HJB: Optimal biomass and quota function when $X_3(0)=0.7$ and $X_3(0)=1.3$.}
\label{simul3hjb3} 
\end{minipage}
\end{figure}

\subsection{Optimisation with a Neural Network with one layer of 15000 Neurons}
The parameters are the same.  The neural network has one layer with 15000 neurons. Conjugate gradient with optimal step size is used and the derivatives are computed with automatic differentiation in reverse mode.

Convergence of the conjugate gradient algorithm is seen on Figure \ref{cost15000}. The norm of the gradient  is reduced from   0.1 to  0.00017. 
\begin{figure}[h!]
\begin{minipage} [b]{0.45\textwidth}
\centering
\begin{tikzpicture}[scale=0.8]
\begin{axis}[legend style={at={(1,1)},anchor=north east}, compat=1.3,xlabel= {$X_0$},ylabel= {average cost}, 
ytick={-0.01,-0.005,0.,0.01,0.02,0.03,0.04,0.05}, grid style=dashed, ymajorgrids=true]
\addplot[thick,solid,color=red,mark=none,mark size=1pt] table [x index=0, y index=1]{results/hjb3b/costs2.txt};
\end{axis}
\end{tikzpicture}
\caption{Values of the average cost versus $X_0$ for $20$ realizations by HJB for the 3 species problem.}
\label{hjb3cost}
\end{minipage}
\hskip 0.5cm
\begin{minipage} [b]{0.45\textwidth}
\centering
\begin{tikzpicture}[scale=0.8]
\begin{axis}[legend style={at={(1,1)},anchor=north east}, compat=1.3,xlabel= {$n$},ylabel= {cost J}]
\addplot[thick,dotted,color=blue,mark=none,mark size=1pt] table [x index=0, y index=1]{fig/convergence.txt};
\addlegendentry{$J^n$}
\end{axis}
\end{tikzpicture}
\caption{\label{cost15000} convergence of the cost function during the optimisation process to solve the 3D problem with a NN with one layer of 150000 neurons.}
\end{minipage}
\end{figure}
Then the results are displayed as above: first the optimal quota vector function, with 9 surfaces on Figures \ref{nn11} to \ref{nn33}.  The fact that the surface $\{t,X_2\}\mapsto u_2(1,X_2,1,t)$ is flat is an obvious indication that the solution found by the NN with one layer is only suboptimal. 

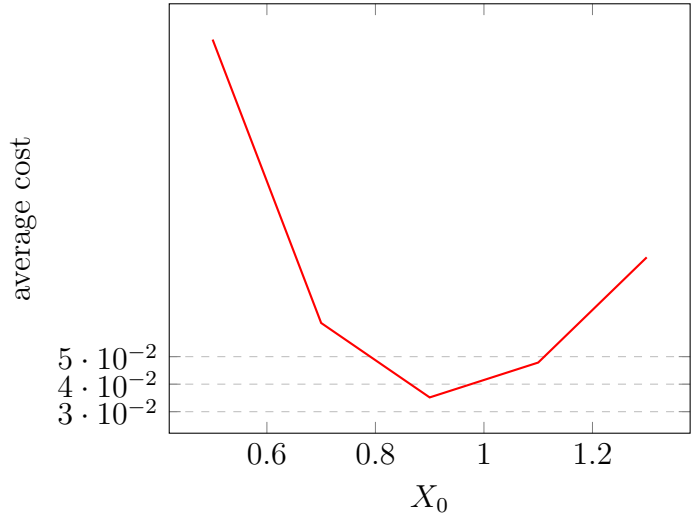
\begin{figure}[h!]
\centering
\begin{tikzpicture}[scale=1]
\begin{axis}[legend style={at={(1,1)},anchor=north east}, compat=1.3,xlabel= {$X_0$},ylabel= {average cost}, 
ytick={-0.01,-0.005,0.,0.01,0.02,0.03,0.04,0.05}, grid style=dashed, ymajorgrids=true]
\addplot[thick,solid,color=red,mark=none,mark size=1pt] table [x index=0, y index=1] {results/D3NN1500Iter200/costs.txt};
\end{axis}
\end{tikzpicture}
\caption{Values of the average cost versus $X_0$ for $10$ realizations of $\d W$ and with the Markovian neural network with 15000 neurons and 10 conjugate gradient iterations.}
\label{costsNN15000}
\end{figure}

Then two types of random trajectories, one starting at $X_0=(0.7,0.7,0.7)^T$, the other at $(1.3,1.3,1.3)^T$.  Results are shown on Figures \ref{nn15000x}, \ref{nn15000y} and \ref{nn15000z}.
The section ends with a display of the cost function versus $\vX_0$ on Figure \ref{costsNN15000}.
The results seem less accurate than with HJB.  Figure \ref{nn22} is surprising! It shows also on Figure  \ref{nn15000y} where the optimal quota does not steer $\vX$ anywhere near $\vX_d$.

\begin{figure}[h!]
\begin{minipage} [b]{0.32\textwidth}
\centering
\begin{tikzpicture}[scale=0.6]
\begin{axis}[legend style={at={(1,1)},anchor=north east}, compat=1.3,xlabel= {$t$},ylabel= {$X$}]
 \addplot3[surf] table [ ] {results/D3NN1500Iter200/surf11.txt};
\addlegendentry{ $u_1(X_1,t)$}
\end{axis}
\end{tikzpicture}
\caption{$X_1,t\mapsto u_1$.}
\label{nn11} 
\end{minipage}\hskip 0.5cm
\begin{minipage} [b]{0.32\textwidth}
\centering
\begin{tikzpicture}[scale=0.6]
\begin{axis}[legend style={at={(1,1)},anchor=north east}, compat=1.3,xlabel= {$t$},ylabel= {$X$}]
 \addplot3[surf] table [ ] {results/D3NN1500Iter200/surf12.txt};
\addlegendentry{ $u_1(X_2,t)$}
\end{axis}
\end{tikzpicture}
\caption{$X_2,t\mapsto u_1$.}
\label{nn12} 
\end{minipage}
\begin{minipage} [b]{0.32\textwidth}
\centering
\begin{tikzpicture}[scale=0.6]
\begin{axis}[legend style={at={(1,1)},anchor=north east}, compat=1.3,xlabel= {$t$},ylabel= {$X$}]
 \addplot3[surf] table [ ] {results/D3NN1500Iter200/surf13.txt};
\addlegendentry{ $u_1(X_3,t,t)$}
\end{axis}
\end{tikzpicture}
\caption{$X_3,t\mapsto u_1$.}
\label{nn13} 
\end{minipage}
\end{figure}

\begin{figure}[h!]
\begin{minipage} [b]{0.3\textwidth}
\centering
\begin{tikzpicture}[scale=0.56]
\begin{axis}[legend style={at={(1,1)},anchor=north east}, compat=1.3,xlabel= {$t$},ylabel= {$X$}]
 \addplot3[surf] table [ ] {results/D3NN1500Iter200/surf21.txt};
\addlegendentry{ $u_2(X_1,1,t)$}
\end{axis}
\end{tikzpicture}
\caption{$X_1,t\mapsto u_2$.}
\label{nn21} 
\end{minipage}\hskip 0.5cm
\begin{minipage} [b]{0.3\textwidth}
\centering
\begin{tikzpicture}[scale=0.56]
\begin{axis}[legend style={at={(1,1)},anchor=north east}, compat=1.3,xlabel= {$t$},ylabel= {$X$}]
 \addplot3[surf] table [ ] {results/D3NN1500Iter200/surf22.txt};
\addlegendentry{ $u_2(X_2,t)$}
\end{axis}
\end{tikzpicture}
\caption{$X_2,t\mapsto u_2$.}
\label{nn22} 
\end{minipage}\hskip 0.5cm
\begin{minipage} [b]{0.3\textwidth}
\centering
\begin{tikzpicture}[scale=0.56]
\begin{axis}[legend style={at={(1,1)},anchor=north east}, compat=1.3,xlabel= {$1$},ylabel= {$X$}]
 \addplot3[surf] table [ ] {results/D3NN1500Iter200/surf23.txt};
\addlegendentry{ $u_2(X_3,t,t)$}
\end{axis}
\end{tikzpicture}
\caption{$X_3,t\mapsto u_2$.}
\label{nn23} 
\end{minipage}
\end{figure}

\begin{figure}[h!]
\begin{minipage} [b]{0.32\textwidth}
\centering
\begin{tikzpicture}[scale=0.6]
\begin{axis}[legend style={at={(1,1)},anchor=north east}, compat=1.3,xlabel= {$t$},ylabel= {$X$}]
 \addplot3[surf] table [ ] {results/D3NN1500Iter200/surf31.txt};
\addlegendentry{ $u_3(X_1,t)$}
\end{axis}
\end{tikzpicture}
\caption{$X_1,t\mapsto u_3$.}
\label{nn31} 
\end{minipage}\hskip 0.5cm
\begin{minipage} [b]{0.32\textwidth}
\centering
\begin{tikzpicture}[scale=0.6]
\begin{axis}[legend style={at={(1,1)},anchor=north east}, compat=1.3,xlabel= {$t$},ylabel= {$X$}]
 \addplot3[surf] table [ ] {results/D3NN1500Iter200/surf32.txt};
\addlegendentry{ $u_3(X_2,t)$}
\end{axis}
\end{tikzpicture}
\caption{$X_2,t\mapsto u_3$.}
\label{nn32} 
\end{minipage}
\begin{minipage} [b]{0.32\textwidth}
\centering
\begin{tikzpicture}[scale=0.6]
\begin{axis}[legend style={at={(1,1)},anchor=north east}, compat=1.3,xlabel= {$t$},ylabel= {$X$}]
 \addplot3[surf] table [ ] {results/D3NN1500Iter200/surf33.txt};
\addlegendentry{ $u_3(X_3,t,t)$}
\end{axis}
\end{tikzpicture}
\caption{$X_3,t\mapsto u_3$.}
\label{nn33} 
\end{minipage}
\end{figure}
\begin{figure}[h!]
\begin{minipage} [b]{0.32\textwidth}
\centering
\begin{tikzpicture}[scale=0.6]
\begin{axis}[legend style={at={(1,1)},anchor=north east}, compat=1.3, ymax=1.4,
  xlabel= {$time$},
  ylabel= {$u_1({X_1}_t,t),~~~~~~{X_1}_t$}]
\addplot[thick,solid,color=black,mark=none,mark size=1pt] table [x index=0, y index=1] {results/D3NN1500Iter200/traj71.txt};
\addlegendentry{ ${u_1}_t$}
\addplot[thin,solid,color=red,mark=none,mark size=1pt] table [x index=0, y index=2] {results/D3NN1500Iter200/traj71.txt};
\addlegendentry{${X_1}_t$}
\addplot[thick,dashed,color=black,mark=none,mark size=1pt] table [x index=0, y index=1] {results/D3NN1500Iter200/traj131.txt};
\addlegendentry{${u_1}_t$}
\addplot[thin,dashed,color=red,mark=none,mark size=1pt] table [x index=0, y index=2] {results/D3NN1500Iter200/traj131.txt};
\addlegendentry{  ${X_1}_t$}
\end{axis}
\end{tikzpicture}
\caption{Optimal biomass and quota function computed with the single layer NN  when $X_0=0.7$ and $X_0=1.3$.}
\label{nn15000x} 
\end{minipage}
\begin{minipage} [b]{0.32\textwidth}
\centering
\begin{tikzpicture}[scale=0.6]
\begin{axis}[legend style={at={(1,1)},anchor=north east}, compat=1.3,
  xlabel= {$time$},
  ylabel= {$u(X_t,t),~~~~~~X_t$}]
\addplot[thick,solid,color=black,mark=none,mark size=1pt] table [x index=0, y index=1] {results/D3NN1500Iter200/traj72.txt};
\addlegendentry{ ${u_2}_t$}
\addplot[thin,solid,color=red,mark=none,mark size=1pt] table [x index=0, y index=2] {results/D3NN1500Iter200/traj72.txt};
\addlegendentry{  ${X_2}_t$}
\addplot[thick,dashed,color=black,mark=none,mark size=1pt] table [x index=0, y index=1] {results/D3NN1500Iter200/traj132.txt};
\addlegendentry{ ${u_2}_t$}
\addplot[thin,dashed,color=red,mark=none,mark size=1pt] table [x index=0, y index=2] {results/D3NN1500Iter200/traj132.txt};
\addlegendentry{  ${X_2}_t$}
\end{axis}
\end{tikzpicture}
\caption{Optimal biomass and quota function computed with the single layer NN  when $X_0=0.7$ and $X_0=1.3$.}
\label{nn15000y} 
\end{minipage}
\begin{minipage} [b]{0.32\textwidth}
\centering
\begin{tikzpicture}[scale=0.6]
\begin{axis}[legend style={at={(1,1)},anchor=north east}, compat=1.3, 
  xlabel= {$time$},
  ylabel= {$u(X_t,t),~~~~~~X_t$}]
\addplot[thick,solid,color=black,mark=none,mark size=1pt] table [x index=0, y index=1] {results/D3NN1500Iter200/traj73.txt};
\addlegendentry{ ${u_3}_t$}
\addplot[thin,solid,color=red,mark=none,mark size=1pt] table [x index=0, y index=2] {results/D3NN1500Iter200/traj73.txt};
\addlegendentry{ ${X_3}_t$}
\addplot[thick,dashed,color=black,mark=none,mark size=1pt] table [x index=0, y index=1] {results/D3NN1500Iter200/traj133.txt};
\addlegendentry{ ${u_3}_t$}
\addplot[thin,dashed,color=red,mark=none,mark size=1pt] table [x index=0, y index=2] {results/D3NN1500Iter200/traj133.txt};
\addlegendentry{ ${X_3}_t$}
\end{axis}
\end{tikzpicture}
\caption{Optimal biomass and quota function computed with the single layer NN  when $X_0=0.7$ and $X_0=1.3$.}
\label{nn15000z} 
\end{minipage}
\end{figure}

\subsection{With two layers}
The same problem is solved with a two layers neural network with 100 neurons on each layer.

 The results are displayed as above: first the optimal quota vector function, with 9 surfaces on Figures \ref{n2n11} to \ref{n2n33}.

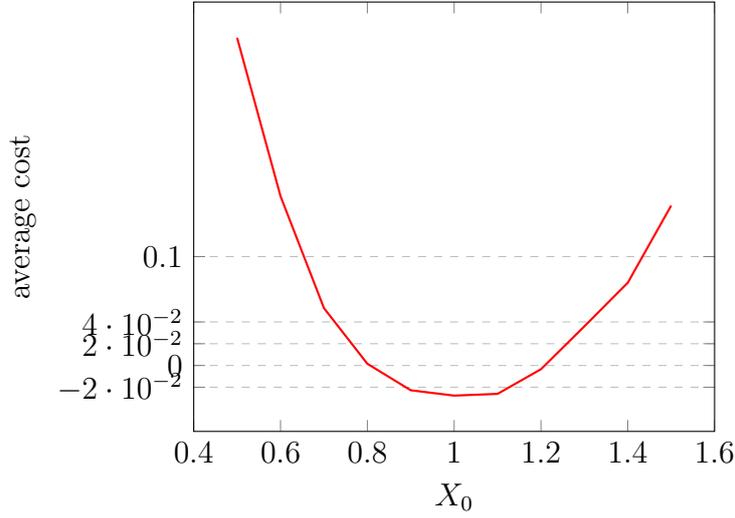
\begin{figure}[h!]
\centering
\begin{tikzpicture}[scale=1]
\begin{axis}[legend style={at={(1,1)},anchor=north east}, compat=1.3,xlabel= {$X_0$},ylabel= {average cost}, 
ytick={-0.02,0,0.02,0.04,0.1}, grid style=dashed, ymajorgrids=true]
\addplot[thick,solid,color=red,mark=none,mark size=1pt] table [x index=0, y index=1] {results/NN2layers/nn_ML_costs.txt};
\end{axis}
\end{tikzpicture}
\caption{Values of the average cost versus $X_0$ for $100$ realizations of $\d W$ and with the Markovian neural network with two layers of neurons.}
\label{costsn2n100}
\end{figure}

Then two types of random trajectories, one starting at $X_0=(0.7,0.7,0.7)^T$, the other at $(1.3,1.3,1.3)^T$.  Results are shown on Figures \ref{n2n100x}, \ref{n2n100y} and \ref{n2n100z}.
The section ends with a display of the cost function versus $\vX_0$ on Figure \ref{costsn2n100}.
\begin{figure}[h!]
\begin{minipage} [b]{0.32\textwidth}
\centering
\begin{tikzpicture}[scale=0.6]
\begin{axis}[legend style={at={(1,1)},anchor=north east}, compat=1.3,xlabel= {$t$},ylabel= {$X$}]
 \addplot3[surf] table [ ] {results/NN2layers/nn_ML_surf11.txt};
\addlegendentry{ $u_1(X_1,t)$}
\end{axis}
\end{tikzpicture}
\caption{$X_1,t\mapsto u_1$.}
\label{n2n11} 
\end{minipage}\hskip 0.5cm
\begin{minipage} [b]{0.32\textwidth}
\centering
\begin{tikzpicture}[scale=0.6]
\begin{axis}[legend style={at={(1,1)},anchor=north east}, compat=1.3,xlabel= {$t$},ylabel= {$X$}]
 \addplot3[surf] table [ ] {results/NN2layers/nn_ML_surf12.txt};
\addlegendentry{ $u_1(X_2,t)$}
\end{axis}
\end{tikzpicture}
\caption{$X_2,t\mapsto u_1$.}
\label{n2n12} 
\end{minipage}
\begin{minipage} [b]{0.32\textwidth}
\centering
\begin{tikzpicture}[scale=0.6]
\begin{axis}[legend style={at={(1,1)},anchor=north east}, compat=1.3,xlabel= {$t$},ylabel= {$X$}]
 \addplot3[surf] table [ ] {results/NN2layers/nn_ML_surf13.txt};
\addlegendentry{ $u_1(X_3,t,t)$}
\end{axis}
\end{tikzpicture}
\caption{$X_3,t\mapsto u_1$.}
\label{n2n13} 
\end{minipage}
\end{figure}

\begin{figure}[h!]
\begin{minipage} [b]{0.3\textwidth}
\centering
\begin{tikzpicture}[scale=0.56]
\begin{axis}[legend style={at={(1,1)},anchor=north east}, compat=1.3,xlabel= {$t$},ylabel= {$X$}]
 \addplot3[surf] table [ ] {results/NN2layers/nn_ML_surf21.txt};
\addlegendentry{ $u_2(X_1,1,t)$}
\end{axis}
\end{tikzpicture}
\caption{$X_1,t\mapsto u_2$.}
\label{n2n21} 
\end{minipage}\hskip 0.5cm
\begin{minipage} [b]{0.3\textwidth}
\centering
\begin{tikzpicture}[scale=0.56]
\begin{axis}[legend style={at={(1,1)},anchor=north east}, compat=1.3,xlabel= {$t$},ylabel= {$X$}]
 \addplot3[surf] table [ ] {results/NN2layers/nn_ML_surf22.txt};
\addlegendentry{ $u_2(X_2,t)$}
\end{axis}
\end{tikzpicture}
\caption{$X_2,t\mapsto u_2$.}
\label{n2n22} 
\end{minipage}\hskip 0.5cm
\begin{minipage} [b]{0.3\textwidth}
\centering
\begin{tikzpicture}[scale=0.56]
\begin{axis}[legend style={at={(1,1)},anchor=north east}, compat=1.3,xlabel= {$1$},ylabel= {$X$}]
 \addplot3[surf] table [ ] {results/NN2layers/nn_ML_surf23.txt};
\addlegendentry{ $u_2(X_3,t,t)$}
\end{axis}
\end{tikzpicture}
\caption{$X_3,t\mapsto u_2$.}
\label{n2n23} 
\end{minipage}
\end{figure}

\begin{figure}[h!]
\begin{minipage} [b]{0.32\textwidth}
\centering
\begin{tikzpicture}[scale=0.6]
\begin{axis}[legend style={at={(1,1)},anchor=north east}, compat=1.3,xlabel= {$t$},ylabel= {$X$}]
 \addplot3[surf] table [ ] {results/NN2layers/nn_ML_surf31.txt};
\addlegendentry{ $u_3(X_1,t)$}
\end{axis}
\end{tikzpicture}
\caption{$X_1,t\mapsto u_3$.}
\label{n2n31} 
\end{minipage}\hskip 0.5cm
\begin{minipage} [b]{0.32\textwidth}
\centering
\begin{tikzpicture}[scale=0.6]
\begin{axis}[legend style={at={(1,1)},anchor=north east}, compat=1.3,xlabel= {$t$},ylabel= {$X$}]
 \addplot3[surf] table [ ] {results/NN2layers/nn_ML_surf32.txt};
\addlegendentry{ $u_3(X_2,t)$}
\end{axis}
\end{tikzpicture}
\caption{$X_2,t\mapsto u_3$.}
\label{n2n32} 
\end{minipage}
\begin{minipage} [b]{0.32\textwidth}
\centering
\begin{tikzpicture}[scale=0.6]
\begin{axis}[legend style={at={(1,1)},anchor=north east}, compat=1.3,xlabel= {$t$},ylabel= {$X$}]
 \addplot3[surf] table [ ] {results/NN2layers/nn_ML_surf33.txt};
\addlegendentry{ $u_3(X_3,t,t)$}
\end{axis}
\end{tikzpicture}
\caption{$X_3,t\mapsto u_3$.}
\label{n2n33} 
\end{minipage}
\end{figure}

\begin{figure}[h!]
\begin{minipage} [b]{0.32\textwidth}
\centering
\begin{tikzpicture}[scale=0.6]
\begin{axis}[legend style={at={(1,1)},anchor=north east}, compat=1.3, ymax=1.4,
  xlabel= {$time$},
  ylabel= {$u_1({X_1}_t,t),~~~~~~{X_1}_t$}]
\addplot[thick,solid,color=black,mark=none,mark size=1pt] table [x index=0, y index=1] {results/NN2layers/nn_ML_traj71.txt};
\addlegendentry{ ${u_1}_t$}
\addplot[thin,solid,color=red,mark=none,mark size=1pt] table [x index=0, y index=2] {results/NN2layers/nn_ML_traj71.txt};
\addlegendentry{${X_1}_t$}
\addplot[thick,dashed,color=black,mark=none,mark size=1pt] table [x index=0, y index=1] {results/NN2layers/nn_ML_traj131.txt};
\addlegendentry{${u_1}_t$}
\addplot[thin,dashed,color=red,mark=none,mark size=1pt] table [x index=0, y index=2] {results/NN2layers/nn_ML_traj131.txt};
\addlegendentry{  ${X_1}_t$}
\end{axis}
\end{tikzpicture}
\caption{Optimal biomass and quota function computed with the two layers NN  when $X_0=0.7$ and $X_0=1.3$.}
\label{n2n100x} 
\end{minipage}
\begin{minipage} [b]{0.32\textwidth}
\centering
\begin{tikzpicture}[scale=0.6]
\begin{axis}[legend style={at={(1,1)},anchor=north east}, compat=1.3,
  xlabel= {$time$},
  ylabel= {$u(X_t,t),~~~~~~X_t$}]
\addplot[thick,solid,color=black,mark=none,mark size=1pt] table [x index=0, y index=1] {results/NN2layers/nn_ML_traj72.txt};
\addlegendentry{ ${u_2}_t$}
\addplot[thin,solid,color=red,mark=none,mark size=1pt] table [x index=0, y index=2] {results/NN2layers/nn_ML_traj72.txt};
\addlegendentry{  ${X_2}_t$}
\addplot[thick,dashed,color=black,mark=none,mark size=1pt] table [x index=0, y index=1] {results/NN2layers/nn_ML_traj132.txt};
\addlegendentry{ ${u_2}_t$}
\addplot[thin,dashed,color=red,mark=none,mark size=1pt] table [x index=0, y index=2] {results/NN2layers/nn_ML_traj132.txt};
\addlegendentry{  ${X_2}_t$}
\end{axis}
\end{tikzpicture}
\caption{Optimal biomass and quota function computed with the two layers NN  when $X_0=0.7$ and $X_0=1.3$.}
\label{n2n100y} 
\end{minipage}
\begin{minipage} [b]{0.32\textwidth}
\centering
\begin{tikzpicture}[scale=0.6]
\begin{axis}[legend style={at={(1,1)},anchor=north east}, compat=1.3, 
  xlabel= {$time$},
  ylabel= {$u(X_t,t),~~~~~~X_t$}]
\addplot[thick,solid,color=black,mark=none,mark size=1pt] table [x index=0, y index=1] {results/NN2layers/nn_ML_traj73.txt};
\addlegendentry{ ${u_3}_t$}
\addplot[thin,solid,color=red,mark=none,mark size=1pt] table [x index=0, y index=2] {results/NN2layers/nn_ML_traj73.txt};
\addlegendentry{ ${X_3}_t$}
\addplot[thick,dashed,color=black,mark=none,mark size=1pt] table [x index=0, y index=1] {results/NN2layers/nn_ML_traj133.txt};
\addlegendentry{ ${u_3}_t$}
\addplot[thin,dashed,color=red,mark=none,mark size=1pt] table [x index=0, y index=2] {results/NN2layers/nn_ML_traj133.txt};
\addlegendentry{ ${X_3}_t$}
\end{axis}
\end{tikzpicture}
\caption{Optimal biomass and quota function computed with the two layers NN  when $X_0=0.7$ and $X_0=1.3$.}
\label{n2n100z} 
\end{minipage}
\end{figure}

Consequently, the neural network with 2 layers give decent results, better than HJB.  The neural network with one layer only gives poor results.

\subsection{Assessment of the results}
Remark \ref{remm1b} allows the construction of a simple 3D solution from a 1D solution.

Let $v,y$ be solution of
\[
\min_v \left\{\int_0^T\E\left[|y-y_d|^2 -\alpha v +\beta |v'|^2\right]\d t ~: ~
\d y_t = y_t(r- y +v +\sigma\d w_t)~,~y(0)=y_0 \right\}
\]
Let $\underline{\bm\kappa}\in\R^{3\times 3}$, $\vu=(v,v,v)^T$, $\vr=(r,r,r)^T$, $\bm\alpha = (\alpha,\alpha,\alpha)^T$, 
$\underline{\bm\sigma}=\underline{\bm\Lambda}[(\sigma,\sigma,\sigma)^T]$, 
$\vY=(y,y,y)^T$ and $\d\vW_t=(\d w_t,\d w_t,\d w_t)^T$.
Then
\[
\d\vY_t= \underline{\bm\Lambda}\left[(\vr-\vu  - \vY)\d t +\underline{\bm \sigma}\d\vW_t \right]\vY
\]
Let $\vX=\underline{\bm\kappa}^{-1}\vY$.
It implies
\[
\d\vX_t = \underline{\bm\kappa}^{-1} \underline{\bm\Lambda}\left[(\vr-\vu  - \underline{\bm\kappa}\vX)\d t +\underline{\bm \sigma}\d\vW_t \right]\underline{\bm\kappa}\vX.
\]
The matrix $\underline{\bm\Lambda}\left[(\vr-\vu  - \underline{\bm\kappa}\vX)\d t +\underline{\bm \sigma}\d\vW_t \right]$ is diagonal  and all terms are equal, so it commutes with $\underline{\bm\kappa}$.  Therefore
\[
\d\vX_t =  \underline{\bm\Lambda}\left[(\vr-\vu  - \underline{\bm\kappa}\vX)\d t +\underline{\bm \sigma}\d\vW_t \right]\vX.
\]
and $\vu$ is solution of
\[
\min_{u\in \U} \bar J:=\E\left[\int_0^T\left[|\underline{\bm\kappa}(\vX-\underline{\bm\kappa}^{-1}\vY_d)|^2\d t  - \bm\alpha \cdot \vu + \beta|\frac{\d\vu}{\d t}|^2\right] \d t \right]
\]
because this is 3 times the cost function of the problem in $y,v$.

Let $\vX,\vu$ be the solution of this problem, then by construction $\vX=\underline{\bm\kappa}^{-1}\vY$ and 
$$
\vu(X_1,X_2,X_3)=(v(Y_1),v(Y_2),v(Y_3))^T=(v(y),v(y),v(y))^T.
$$
For instance, $\vu_1(1,X_2,1)=v((\underline{\bm\kappa}(1,X_2,1)^T)_1)$.

Now we build on the fact that $v$ is approximatively bang-bang and equal to 0.5 when $y<1$ and 1 otherwise.

Example
\[
\underline{\bm\kappa}=
\left(\begin{matrix}
1.2&-0.1&0\cr 0.2&1.2&0\cr 0.1&0.1&1\cr
\end{matrix}\right)
\quad
\underline{\bm\kappa}^{-1}=
\left(\begin{matrix}
0.822&0.0685& 0\cr
-0.137& 0.822& 0\cr
-0.0685&-0.089&1
\end{matrix}\right)
\quad
\underline{\bm\kappa}^{-1}\left(\begin{matrix}
1\cr 1\cr 1\cr
\end{matrix}\right)
=\left(\begin{matrix}
1.16\cr 0.89\cr 1.09\cr
\end{matrix}\right)\times 0.77
\]
The scaling $0.77$ plays no role. 
\[
\vu(X_1,X_{d2},X_{d3})=
v\left( \underline{\bm\kappa}\left(\begin{matrix}
X_1\cr 0.89 \times 0.77 = 0.685 \cr 1.09 \times 0.77 = 0.839 \cr
\end{matrix}\right) \right)
= v\left(\left(\begin{matrix}
1.2{X_1} - 0.69\cr
0.2 X_1 + 0.82 \cr
0.1 X_1 +0.91\cr
\end{matrix}\right) \right)
\]
Hence $\vu_1(X_1,X_{d2},X_{d3})$ is expected to be bang-bang at $X_1=0.56$, $\vu_2(X_1,1,1)$    is expected to be bang bang at $X_1 = 0.9$, and $\vu_3(X_1,1,1)$  is expected to be bang bang at $X_1 = 0.9$.
\[
\vu(X_{d1},X_2,X_{d3})=
v\left(\underline{\bm\kappa}\left(\begin{matrix}
1.16 \times 0.77 = 0.89 \cr X_2\cr 1.09 \times 0.77 = 0.89 \cr
\end{matrix}\right) \right)
= v\left(\left(\begin{matrix}
1.07-0.1 X_2\cr
0.18+1.2 X_2\cr
0.98 + 0.1 X_2 \cr
\end{matrix}\right) \right)
\]
Hence $\vu_1(X_{d1},X_2,X_{d3})$ is expected to be bang bang at $X_2=0.7$, and $\vu_2(X_{d1},X_2,X_{d3})$ is expected to be bang-bang at $X_2=0.68$ and $\vu_3(X_{d1},X_2,X_{d3})$ is expected equal to be bang bang at $X_2 = 0.2$. 

\[
\vu(X_{d1},X_{d2},X_3)
=
	v\left(\underline{\bm\kappa}\left(\begin{matrix}
	1.16 \times 0.77 = 0.89 \cr 0.89 \times 0.77 = 0.685 \cr X_3\cr
	\end{matrix}\right) \right)
=
	v\left(\left(\begin{matrix}
	1.0\cr
	1.0\cr
	0.16+X_3 \cr
	\end{matrix}\right)\right)
\]
Hence $\vu_1(X_{d1},X_{d2},X_3)$  and $\vu_2(X_{d1},X_{d2},X_3)$, are expected equal to $1$   everywhere and $\vu_3(X_{d1},X_{d2},X_3))$ is expected to be bang-bang at $X_3=0.84$.

\subsection{Optimisation with HJB when $X_d=\underline{\bm\kappa}^{-1}{\bf 1}$}
The same numerical test was done with HJB but with $\vX_d=(1.16,0.089,1.09)^T$.
The results are shown on Figures \ref{Xdnn11} to \ref{Xdnn15000z}.  The lowest cost is $0.04$, for $X_0=1$.

The solution is nearer to the one constructed above but still fairly different. The Neural Networks performances are poor on this test. It is likely that they produced a local minimum.  

\begin{figure}[h!]
\begin{minipage} [b]{0.32\textwidth}
\centering
\begin{tikzpicture}[scale=0.6]
\begin{axis}[legend style={at={(1,1)},anchor=north east}, compat=1.3,xlabel= {$t$},ylabel= {$X$}]
 \addplot3[surf] table [ ] {results/nnD3Xdnot1/surf11.txt};
\addlegendentry{ $u_1(X_1,t)$}
\end{axis}
\end{tikzpicture}
\caption{3 species: $X_1,t\mapsto u_1$.}
\label{Xdnn11} 
\end{minipage}\hskip 0.5cm
\begin{minipage} [b]{0.32\textwidth}
\centering
\begin{tikzpicture}[scale=0.6]
\begin{axis}[legend style={at={(1,1)},anchor=north east}, compat=1.3,xlabel= {$t$},ylabel= {$X$}]
 \addplot3[surf] table [ ] {results/nnD3Xdnot1/surf12.txt};
\addlegendentry{ $u_1(X_2,t)$}
\end{axis}
\end{tikzpicture}
\caption{3 species: $X_2,t\mapsto u_1$.}
\label{Xdnn12} 
\end{minipage}
\begin{minipage} [b]{0.32\textwidth}
\centering
\begin{tikzpicture}[scale=0.6]
\begin{axis}[legend style={at={(1,1)},anchor=north east}, compat=1.3,xlabel= {$t$},ylabel= {$X$}]
 \addplot3[surf] table [ ] {results/nnD3Xdnot1/surf13.txt};
\addlegendentry{ $u_1(X_3,t)$}
\end{axis}
\end{tikzpicture}
\caption{3 species: $X_3,t\mapsto u_1$.}
\label{Xdnn13} 
\end{minipage}
\end{figure}

\begin{figure}[h!]
\begin{minipage} [b]{0.3\textwidth}
\centering
\begin{tikzpicture}[scale=0.56]
\begin{axis}[legend style={at={(1,1)},anchor=north east}, compat=1.3,xlabel= {$t$},ylabel= {$X$}]
 \addplot3[surf] table [ ] {results/nnD3Xdnot1/surf21.txt};
\addlegendentry{ $u_2(X_1,t)$}
\end{axis}
\end{tikzpicture}
\caption{3 species: $X_1,t\mapsto u_2$.}
\label{Xdnn21} 
\end{minipage}\hskip 0.5cm
\begin{minipage} [b]{0.3\textwidth}
\centering
\begin{tikzpicture}[scale=0.56]
\begin{axis}[legend style={at={(1,1)},anchor=north east}, compat=1.3,xlabel= {$t$},ylabel= {$X$}]
 \addplot3[surf] table [ ] {results/nnD3Xdnot1/surf22.txt};
\addlegendentry{ $u_2(X_2,t)$}
\end{axis}
\end{tikzpicture}
\caption{3 species: $X_2,t\mapsto u_2$.}
\label{Xdnn22} 
\end{minipage}\hskip 0.5cm
\begin{minipage} [b]{0.3\textwidth}
\centering
\begin{tikzpicture}[scale=0.56]
\begin{axis}[legend style={at={(1,1)},anchor=north east}, compat=1.3,xlabel= {$1$},ylabel= {$X$}]
 \addplot3[surf] table [ ] {results/nnD3Xdnot1/surf23.txt};
\addlegendentry{ $u_2(X_3,t)$}
\end{axis}
\end{tikzpicture}
\caption{3 species: $X_3,t\mapsto u_2$.}
\label{Xdnn23} 
\end{minipage}
\end{figure}

\begin{figure}[h!]
\begin{minipage} [b]{0.32\textwidth}
\centering
\begin{tikzpicture}[scale=0.6]
\begin{axis}[legend style={at={(1,1)},anchor=north east}, compat=1.3,xlabel= {$t$},ylabel= {$X$}]
 \addplot3[surf] table [ ] {results/nnD3Xdnot1/surf31.txt};
\addlegendentry{ $u_3(X_1,t)$}
\end{axis}
\end{tikzpicture}
\caption{3 species: $X_1,t\mapsto u_3$.}
\label{Xdnn31} 
\end{minipage}\hskip 0.5cm
\begin{minipage} [b]{0.32\textwidth}
\centering
\begin{tikzpicture}[scale=0.6]
\begin{axis}[legend style={at={(1,1)},anchor=north east}, compat=1.3,xlabel= {$t$},ylabel= {$X$}]
 \addplot3[surf] table [ ] {results/nnD3Xdnot1/surf32.txt};
\addlegendentry{ $u_3(X_2,t)$}
\end{axis}
\end{tikzpicture}
\caption{3 species:: $X_2,t\mapsto u_3$.}
\label{Xdnn32} 
\end{minipage}
\begin{minipage} [b]{0.32\textwidth}
\centering
\begin{tikzpicture}[scale=0.6]
\begin{axis}[legend style={at={(1,1)},anchor=north east}, compat=1.3,xlabel= {$t$},ylabel= {$X$}]
 \addplot3[surf] table [ ] {results/nnD3Xdnot1/surf33.txt};
\addlegendentry{ $u_3(X_3,t)$}
\end{axis}
\end{tikzpicture}
\caption{3 species: $X_3,t\mapsto u_3$.}
\label{Xdnn33} 
\end{minipage}
\end{figure}
\begin{figure}[h!]
\begin{minipage} [b]{0.32\textwidth}
\centering
\begin{tikzpicture}[scale=0.6]
\begin{axis}[legend style={at={(1,1)},anchor=north east}, compat=1.3, ymax=1.4,
  xlabel= {$time$},
  ylabel= {$u_1({X_1}_t,t),~~~~~~{X_1}_t$}]
\addplot[thick,solid,color=black,mark=none,mark size=1pt] table [x index=0, y index=1] {results/nnD3Xdnot1/traj71.txt};
\addlegendentry{ ${u_1}_t$}
\addplot[thin,solid,color=red,mark=none,mark size=1pt] table [x index=0, y index=2] {results/nnD3Xdnot1/traj71.txt};
\addlegendentry{${X_1}_t$}
\addplot[thick,dashed,color=black,mark=none,mark size=1pt] table [x index=0, y index=1] {results/nnD3Xdnot1/traj131.txt};
\addlegendentry{${u_1}_t$}
\addplot[thin,dashed,color=red,mark=none,mark size=1pt] table [x index=0, y index=2] {results/nnD3Xdnot1/traj131.txt};
\addlegendentry{  ${X_1}_t$}
\end{axis}
\end{tikzpicture}
\caption{3 species: Optimal biomass and quota function computed with the single layer NN  when $X_0=0.7$ and $X_0=1.3$.}
\label{Xdnn15000x} 
\end{minipage}
\begin{minipage} [b]{0.32\textwidth}
\centering
\begin{tikzpicture}[scale=0.6]
\begin{axis}[legend style={at={(1,1)},anchor=north east}, compat=1.3,
  xlabel= {$time$},
  ylabel= {$u(X_t,t),~~~~~~X_t$}]
\addplot[thick,solid,color=black,mark=none,mark size=1pt] table [x index=0, y index=1] {results/nnD3Xdnot1/traj72.txt};
\addlegendentry{ ${u_2}_t$}
\addplot[thin,solid,color=red,mark=none,mark size=1pt] table [x index=0, y index=2] {results/nnD3Xdnot1/traj72.txt};
\addlegendentry{  ${X_2}_t$}
\addplot[thick,dashed,color=black,mark=none,mark size=1pt] table [x index=0, y index=1] {results/nnD3Xdnot1/traj132.txt};
\addlegendentry{ ${u_2}_t$}
\addplot[thin,dashed,color=red,mark=none,mark size=1pt] table [x index=0, y index=2] {results/nnD3Xdnot1/traj132.txt};
\addlegendentry{  ${X_2}_t$}
\end{axis}
\end{tikzpicture}
\caption{3 species:Optimal biomass and quota function computed with the single layer NN  when $X_0=0.7$ and $X_0=1.3$.}
\label{Xdnn15000y} 
\end{minipage}
\begin{minipage} [b]{0.32\textwidth}
\centering
\begin{tikzpicture}[scale=0.6]
\begin{axis}[legend style={at={(1,1)},anchor=north east}, compat=1.3, 
  xlabel= {$time$},
  ylabel= {$u(X_t,t),~~~~~~X_t$}]
\addplot[thick,solid,color=black,mark=none,mark size=1pt] table [x index=0, y index=1] {results/nnD3Xdnot1/traj73.txt};
\addlegendentry{ ${u_3}_t$}
\addplot[thin,solid,color=red,mark=none,mark size=1pt] table [x index=0, y index=2] {results/nnD3Xdnot1/traj73.txt};
\addlegendentry{ ${X_3}_t$}
\addplot[thick,dashed,color=black,mark=none,mark size=1pt] table [x index=0, y index=1] {results/nnD3Xdnot1/traj133.txt};
\addlegendentry{ ${u_3}_t$}
\addplot[thin,dashed,color=red,mark=none,mark size=1pt] table [x index=0, y index=2] {results/nnD3Xdnot1/traj133.txt};
\addlegendentry{ ${X_3}_t$}
\end{axis}
\end{tikzpicture}
\caption{3 species:Optimal biomass and quota function computed with the single layer NN  when $X_0=0.7$ and $X_0=1.3$.}
\label{Xdnn15000z} 
\end{minipage}
\end{figure}

\section{Numerical results: 5 species}
The great advantage of Neural Network optimisation is that it scales well with dimensions.  So to show that it is possible, with the same computer code with very few modifications we computed with the one-layer NN with 150000 neurons a case with 5 species. 20 iterations of conjugate gradients were done.

All parameters are as above except the species correlation matrix:

\[\underline{\bm\kappa}=\left(\begin{matrix}
1.2&-0.1&0.0&0.0&-0.1\cr
0.2&1.2&0.0&0.0&-0.1\cr
0.0&0.2&1.2&-0.1&0.0\cr
0.0&0.0&0.1&1.2&0.0\cr
0.1&0.1&0.0&0.0&1.2\cr
\end{matrix}\right)
\]

Only components 1,2,5 are shown see Figures \ref{D5nn11} to \ref{D5nn33}. On Figures \ref{D5nn15000x}, \ref{D5nn15000y} and \ref{D5nn15000z} trajectories with  optimal quotas show that the NN solution is in general driving $\vX$ towards $\vX_d$. 
\begin{figure}[h!]
\begin{minipage} [b]{0.32\textwidth}
\centering
\begin{tikzpicture}[scale=0.6]
\begin{axis}[legend style={at={(1,1)},anchor=north east}, compat=1.3,xlabel= {$t$},ylabel= {$X$}]
 \addplot3[surf] table [ ] {results/D5one20/surf11.txt};
\addlegendentry{ $u_1(X_1,t)$}
\end{axis}
\end{tikzpicture}
\caption{5 species: $X_1,t\mapsto u_1$.}
\label{D5nn11} 
\end{minipage}\hskip 0.5cm
\begin{minipage} [b]{0.32\textwidth}
\centering
\begin{tikzpicture}[scale=0.6]
\begin{axis}[legend style={at={(1,1)},anchor=north east}, compat=1.3,xlabel= {$t$},ylabel= {$X$}]
 \addplot3[surf] table [ ] {results/D5one20/surf12.txt};
\addlegendentry{ $u_1(X_2,t)$}
\end{axis}
\end{tikzpicture}
\caption{5 species: $X_2,t\mapsto u_1$.}
\label{D5nn12} 
\end{minipage}
\begin{minipage} [b]{0.32\textwidth}
\centering
\begin{tikzpicture}[scale=0.6]
\begin{axis}[legend style={at={(1,1)},anchor=north east}, compat=1.3,xlabel= {$t$},ylabel= {$X$}]
 \addplot3[surf] table [ ] {results/D5one20/surf13.txt};
\addlegendentry{ $u_1(X_D,t)$}
\end{axis}
\end{tikzpicture}
\caption{5 species: $X_D,t\mapsto u_1$.}
\label{D5nn13} 
\end{minipage}
\end{figure}

\begin{figure}[h!]
\begin{minipage} [b]{0.3\textwidth}
\centering
\begin{tikzpicture}[scale=0.56]
\begin{axis}[legend style={at={(1,1)},anchor=north east}, compat=1.3,xlabel= {$t$},ylabel= {$X$}]
 \addplot3[surf] table [ ] {results/D5one20/surf21.txt};
\addlegendentry{ $u_2(X_1,t)$}
\end{axis}
\end{tikzpicture}
\caption{5 species: $X_1,t\mapsto u_2$.}
\label{D5nn21} 
\end{minipage}\hskip 0.5cm
\begin{minipage} [b]{0.3\textwidth}
\centering
\begin{tikzpicture}[scale=0.56]
\begin{axis}[legend style={at={(1,1)},anchor=north east}, compat=1.3,xlabel= {$t$},ylabel= {$X$}]
 \addplot3[surf] table [ ] {results/D5one20/surf22.txt};
\addlegendentry{ $u_2(X_2,t)$}
\end{axis}
\end{tikzpicture}
\caption{5 species: $X_2,t\mapsto u_2$.}
\label{D5nn22} 
\end{minipage}\hskip 0.5cm
\begin{minipage} [b]{0.3\textwidth}
\centering
\begin{tikzpicture}[scale=0.56]
\begin{axis}[legend style={at={(1,1)},anchor=north east}, compat=1.3,xlabel= {$1$},ylabel= {$X$}]
 \addplot3[surf] table [ ] {results/D5one20/surf23.txt};
\addlegendentry{ $u_2(X_3,t)$}
\end{axis}
\end{tikzpicture}
\caption{5 species: $X_3,t\mapsto u_2$.}
\label{D5nn23} 
\end{minipage}
\end{figure}

\begin{figure}[h!]
\begin{minipage} [b]{0.32\textwidth}
\centering
\begin{tikzpicture}[scale=0.6]
\begin{axis}[legend style={at={(1,1)},anchor=north east}, compat=1.3,xlabel= {$t$},ylabel= {$X$}]
 \addplot3[surf] table [ ] {results/D5one20/surf31.txt};
\addlegendentry{ $u_3(X_1,t)$}
\end{axis}
\end{tikzpicture}
\caption{5 species: $X_1,t\mapsto u_3$.}
\label{D5nn31} 
\end{minipage}\hskip 0.5cm
\begin{minipage} [b]{0.32\textwidth}
\centering
\begin{tikzpicture}[scale=0.6]
\begin{axis}[legend style={at={(1,1)},anchor=north east}, compat=1.3,xlabel= {$t$},ylabel= {$X$}]
 \addplot3[surf] table [ ] {results/D5one20/surf32.txt};
\addlegendentry{ $u_D(X_2,t)$}
\end{axis}
\end{tikzpicture}
\caption{5 species: $X_2,t\mapsto u_3$.}
\label{D5nn32} 
\end{minipage}
\begin{minipage} [b]{0.32\textwidth}
\centering
\begin{tikzpicture}[scale=0.6]
\begin{axis}[legend style={at={(1,1)},anchor=north east}, compat=1.3,xlabel= {$t$},ylabel= {$X$}]
 \addplot3[surf] table [ ] {results/D5one20/surf33.txt};
\addlegendentry{ $u_D(X_D,t)$}
\end{axis}
\end{tikzpicture}
\caption{5 species: $X_3,t\mapsto u_3$.}
\label{D5nn33} 
\end{minipage}
\end{figure}
\begin{figure}[h!]
\begin{minipage} [b]{0.32\textwidth}
\centering
\begin{tikzpicture}[scale=0.6]
\begin{axis}[legend style={at={(1,1)},anchor=north east}, compat=1.3, ymax=1.4,
  xlabel= {$time$},
  ylabel= {$u_1({X_1}_t,t),~~~~~~{X_1}_t$}]
\addplot[thick,solid,color=black,mark=none,mark size=1pt] table [x index=0, y index=1] {results/D5one20/traj71.txt};
\addlegendentry{ ${u_1}_t$}
\addplot[thin,solid,color=red,mark=none,mark size=1pt] table [x index=0, y index=2] {results/D5one20/traj71.txt};
\addlegendentry{${X_1}_t$}
\addplot[thick,dashed,color=black,mark=none,mark size=1pt] table [x index=0, y index=1] {results/D5one20/traj131.txt};
\addlegendentry{${u_1}_t$}
\addplot[thin,dashed,color=red,mark=none,mark size=1pt] table [x index=0, y index=2] {results/D5one20/traj131.txt};
\addlegendentry{  ${X_1}_t$}
\end{axis}
\end{tikzpicture}
\caption{5 species: Optimal biomass and quota function computed with the single layer NN  when $X_0=0.7$ and $X_0=1.3$.}
\label{D5nn15000x} 
\end{minipage}
\begin{minipage} [b]{0.32\textwidth}
\centering
\begin{tikzpicture}[scale=0.6]
\begin{axis}[legend style={at={(1,1)},anchor=north east}, compat=1.3,
  xlabel= {$time$},
  ylabel= {$u(X_t,t),~~~~~~X_t$}]
\addplot[thick,solid,color=black,mark=none,mark size=1pt] table [x index=0, y index=1] {results/D5one20/traj72.txt};
\addlegendentry{ ${u_2}_t$}
\addplot[thin,solid,color=red,mark=none,mark size=1pt] table [x index=0, y index=2] {results/D5one20/traj72.txt};
\addlegendentry{  ${X_2}_t$}
\addplot[thick,dashed,color=black,mark=none,mark size=1pt] table [x index=0, y index=1] {results/D5one20/traj132.txt};
\addlegendentry{ ${u_2}_t$}
\addplot[thin,dashed,color=red,mark=none,mark size=1pt] table [x index=0, y index=2] {results/D5one20/traj132.txt};
\addlegendentry{  ${X_2}_t$}
\end{axis}
\end{tikzpicture}
\caption{5 species: Optimal biomass and quota function computed with the single layer NN  when $X_0=0.7$ and $X_0=1.3$.}
\label{D5nn15000y} 
\end{minipage}
\begin{minipage} [b]{0.32\textwidth}
\centering
\begin{tikzpicture}[scale=0.6]
\begin{axis}[legend style={at={(1,1)},anchor=north east}, compat=1.3, 
  xlabel= {$time$},
  ylabel= {$u(X_t,t),~~~~~~X_t$}]
\addplot[thick,solid,color=black,mark=none,mark size=1pt] table [x index=0, y index=1] {results/D5one20/traj73.txt};
\addlegendentry{ ${u_D}_t$}
\addplot[thin,solid,color=red,mark=none,mark size=1pt] table [x index=0, y index=2] {results/D5one20/traj73.txt};
\addlegendentry{ ${X_D}_t$}
\addplot[thick,dashed,color=black,mark=none,mark size=1pt] table [x index=0, y index=1] {results/D5one20/traj133.txt};
\addlegendentry{ ${u_D}_t$}
\addplot[thin,dashed,color=red,mark=none,mark size=1pt] table [x index=0, y index=2] {results/D5one20/traj133.txt};
\addlegendentry{ ${X_D}_t$}
\end{axis}
\end{tikzpicture}
\caption{5 species: Optimal biomass and quota function computed with the single layer NN  when $X_0=0.7$ and $X_0=1.3$.}
\label{D5nn15000z} 
\end{minipage}
\end{figure}
The minimum of the cost function is 0.038 for $\vX_0=1$.

\section{Conclusion}
Control of the biomass of a fishing site has been here a mathematical opportunity to test the numerical methods at hand. One of the advantage of the model is that it is meaningful in any dimension, the number of species. Thus it is a testbed for stochastic control numerical methods.

We have compared Hamilton-Jacobi-Bellman solutions and stochastic dynamic programming  with two implementation of a neural network based optimisation.

The later is conceptually very simple and the computer libraries of AI and Automatic Differentiation can be used.  But they need to be validated and it is the object of this article.

Stochastic dynamic programming is difficult to use numerically beyond dimension 2 and HJB solutions do not scale beyond dimension 4, unless sophisticated discretisation tools are used like sparse grids.  

Neural network based optimisation can be used for large dimension problems but assessing the precision of the answer seems difficult.  Here in dimension 3 the one layer network did not work well and with the two layers network we could observe discrepancies with HJB solutions. In dimension 5, the numerical solution seems reasonable but it is probably suboptimal.

\end{document}